\documentclass[article]{IEEEtran}
\usepackage{color}
\usepackage[pdftex]{graphicx}
\graphicspath{{fig/}{jpeg/}}
\usepackage[cmex10]{amsmath}
\usepackage{amssymb}
\usepackage{epstopdf}
\usepackage{mathtools}
\usepackage{amssymb,bbm}
\usepackage{algorithm,algorithmic,multirow,hhline}
\usepackage{amssymb,lipsum}
\let\labelindent\relax
\usepackage{hyperref,enumitem}
\usepackage{comment}
\usepackage[font=small]{caption}
\usepackage{subcaption} 
\usepackage{mathtools}
\input{mysymbol.sty}
\usepackage{needspace}

\usepackage{theorem}

\newcommand{\y}{{\bf y}}
\newcommand{\x}{{\bf x}}

\newtheorem{theorem}{Theorem}

\newtheorem{lemma}{Lemma}

\newtheorem{assumption}{Assumption}

\date{}
\def\@fnsymbol#1{\ensuremath{\ifcase#1\or *\or \dagger\or \ddagger\or
   \mathsection\or \mathparagraph\or \|\or **\or \dagger\dagger
   \or \ddagger\ddagger \else\@ctrerr\fi}}
\def \QED {$\blacksquare$}

 \newcommand{\INDSTATE}[1][1]{\STATE\hspace{3mm}}
\newcommand{\INDSTATED}[1][1]{\STATE\hspace{6mm}}

\title{Nearly Consistent Finite Particle Estimates in Importance Sampling}

\begin{document}
\title{Nearly Consistent Finite Particle Estimates\\ in Streaming Importance Sampling}

%\titlerunning{Short form of title}        % if too long for running head
%
%\author{Alec Koppel$^\star$         \and
%        Amrit Singh Bedi$^\star$ \and
%        V\'{i}ctor Elvira \and
%        Brian M. Sadler %etc.
%}

\author{Alec Koppel\IEEEauthorrefmark{1},
Amrit Singh Bedi\IEEEauthorrefmark{1},
 Brian M. Sadler\IEEEauthorrefmark{1}, and
  V\'{i}ctor Elvira\IEEEauthorrefmark{2}
\thanks{ \noindent  A. Koppel and A.S. Bedi contributed equally to this work. A preliminary version of this work has appeared as \cite{bedi2019compressed,koppel2019compressed} with no proofs and a different (heuristic) projection rule for quantifying the difference between distributions. 

\noindent \IEEEauthorrefmark{1}U.S. Army Research Laboratory, Adelphi, MD 20783, USA. E-mails: \{alec.e.koppel.civ, brian.m.sadler6.civ\}@mail.mil, amrit0714@gmail.com

\noindent \IEEEauthorrefmark{2}School of Mathematics, University of Edinburgh, Peter Guthrie Tait Road, Edinburgh, EH9 3FD
United Kingdom:  victor.elvira@ed.ac.uk}
\vspace{-0cm}}

\maketitle

\begin{abstract}
In Bayesian inference, we seek to compute information about random variables such as moments or quantiles on the basis of {available data} and prior information. When the distribution of random variables is {intractable}, Monte Carlo (MC) sampling is usually required. {Importance sampling is a standard MC tool that approximates this unavailable distribution with a set of weighted samples.} This procedure is asymptotically consistent as the number of MC samples (particles) go to infinity. However, retaining infinitely many particles is intractable. Thus, we propose a way to only  keep a \emph{finite representative subset} of particles and their augmented importance weights that is \emph{nearly consistent}. To do so in {an online manner}, we (1) embed the posterior density estimate in a reproducing kernel Hilbert space (RKHS) through its kernel mean embedding;  and (2) sequentially project this RKHS element onto a lower-dimensional subspace in RKHS using the maximum mean discrepancy, an integral probability metric. Theoretically, we establish that this scheme results in a bias determined by a compression parameter, which yields a tunable tradeoff between consistency and memory. In experiments,  we observe the compressed estimates achieve comparable performance to the dense ones with substantial reductions in representational complexity.
\end{abstract}

%%%%%%%%%%%%%%%%%%%%%%%%%%%%%%%%%%%%%%%%%%%%%%%%%%%%%%%%%%%%%%%%%%%%%%%%%%%%%%%%%%%%%%%%%%%%%%%%%%%%%%%%%%%%%%%%%%%%%%%%%%%%%%%%%%%%%%%%%%%%%%%%%%%%%%%S  E  C  T  I  O  N %%%%%%%%%%%%%%%%%%%%%%%%%%%%%%%%%%%%%%%%%%%%%%%%%%%%%%%%%%%%%%%%%%%%%%%%%%%%%%%%%%%%%%%%%%%%%%%%%%%%%%%%
%%%%%%%%%%%%%%%%%%%%%%%%%%%%%%%%%%%%%%%%%%%%%%%%%%%%%%%%%%%%%%%%%%
\section{Introduction}\label{sec:introduction}
Bayesian inference is devoted to estimating unknowns by considering as them random variables and constructing a posterior distribution. This posterior distribution incorporates the information of available observations (likelihood function) which is merged with prior knowledge about the  unknown (prior distribution) \cite{robert2013monte}. % are devoted to  is a method where one estimates an unknown distribution on the basis of a prior and observations from a likelihood model .} 
Its application is widespread, spanning statistics \cite{robert2007bayesian}, signal processing \cite{candy2016bayesian}, machine learning \cite{murphy2012machine}, genetics \cite{beaumont2004bayesian}, communications \cite{7227082}, econometrics \cite{jangmin2004stock}, robotics \cite{wurm2010octomap}, among many other examples.
{Bayesian inference is only possible in a very small subset of problems (e.g., when the underlying model between observations and the hidden state is linear and corrupted by Gaussian noise  \cite{kalman}). In some cases, it is possible to linearize nonlinear models and still obtain closed-form (but approximate) solutions  \cite{hoshiya1984structural}. Beyond linearity, Bayesian inference may be tackled by Gaussian Processes when smoothness and unimodality are present \cite{rasmussen2004gaussian}.
However, in most models of interest, closeed-form expressions are not available, and the posterior distribution of the unknowns must be approximated.}
%Unfortunately, when the distributional models are both nonlinear and multi-modal, the aforementioned approaches break down.

The standard approximation methodologies in Bayesian inference are either (a) Monte Carlo (MC) algorithms \cite{neal1993probabilistic}, which include Markov chain Monte Carlo (MCMC) methods \cite{martino2014metropolis}, importance sampling (IS) \cite{elvira2019generalized}, and particle filters (PFs) \cite{djuric2003particle}; or (b) variational algorithms \cite{blei2017variational}. The later approach approximates the posterior by using an optimization algorithm to select within a parametrized family (e.g., by minimizing the KL divergence). Variational algorithms can optimize the parameters sequentially via the solution of a stochastic optimization problem \cite{minka2013expectation}. In recent years, efforts to approximate the expectation by sampling have been investigated, e.g., see \cite{hoffman2013stochastic}.
Unfortunately, unless the prior belongs to a simple exponential family, the parametric update is defined by a non-convex objective, meaning that asymptotic unbiasedness is mostly beyond reach. Mixture models have been considered \cite{liu2016stein}, but their convergence is challenging to characterize and the subject of recent work on resampling  \cite{zhang2018stochastic}.

In contrast, MC is a general approach to Bayesian inference based upon sampling (i.e., the simulation of samples/particles), and is known to generate (weighted) samples that eventually converge to the true distribution/quantity of interest \cite{agapiou2017importance,second}. However, the scalability of Monte Carlo methods is still an ongoing challenge from several angles, in that to obtain consistency, the number of particles must go to infinity, and typically its scaling is exponential in the parameter dimension \cite{bengtsson2008curse}. These scalability problems are the focus of this work, which we specifically study in the context of importance sampling (IS), a very relevant family of MC methods.
We focus on IS, as compared with Markov chain Monte Carlo (MCMC), due to advantages such as, e.g., no burn-in period, simple parallelization \cite{robert2018accelerating}, built-in approximation of the normalizing constant that is  useful in many practical problems (e.g., model selection), and the ability to incorporate adaptive mechanisms without compromising its convergence \cite{bugallo2017adaptive}.  IS methods approximate expectations of arbitrary functions of the unknown parameter via weighted samples generated from one or several proposal densities \cite{first,elvira2019generalized}. Their convergence in terms of the integral approximation error, which vanishes  as the number of samples increases, has been a topic of interest recently \cite{first,agapiou2017importance}, and various statistics to quantify their performance have been proposed \cite{elvira2018rethinking}. 

Our goal is to alleviate the dependence of the convergence rate on the representational complexity of the posterior estimate. To do so, we propose projecting the posterior density onto a finite statistically significant subset of particles {after every particle is generated.\footnote{Doing so may be generalized to scenarios where projection can be performed after generating a fixed number $T$ of particles, a form of mini-batching, but this is omitted for simplicity.}} However, doing so directly in a measure space is often intractable to evaluate. By embedding this density in a reproducing kernel Hilbert space (RKHS) via its kernel mean embedding  \cite{sriperumbudur2010hilbert}, we may  compute projections of distributions via parametric computations involving the RKHS. More specifically, kernel mean embeddings extend the idea of feature mapping to spaces of probability distributions, which, under some regularity conditions \cite[Sec. 3.8]{sriperumbudur2010hilbert}, admits a bijection between probability distributions and RKHS elements.   

{\noindent \bf Contributions.} Based upon this insight, we propose a compression scheme that operates online within importance sampling, sequentially deciding which particles are statistically significant for the integral estimation. To do so, we invoke the idea of distribution embedding \cite{sriperumbudur2010hilbert} and map our unnormalized distributional estimates to RKHS, {in contrast to \cite{martino2021compressed,martino2021astronomy}}. We show that the empirical kernel mean embedding estimates in RKHS are parameterized by the importance weights and particles. Then, we propose to sequentially project embedding estimates onto subspaces of the dictionary of particles, where the dictionary is greedily selected to ensure compressed estimates are close to the uncompressed one (according to some metric). This greedy selection is achieved with a custom variant of matching pursuit \cite{Pati1993} based upon the Maximum Mean Discrepancy, which is an easily computable way to evaluate an integral probability metric by virtue of the RKHS mean embedding. The underpinning of this idea is similar to gradient projections in optimization, which has been exploited recently to surmount memory challenges in kernel and Gaussian process regression \cite{koppelconsistent,Koppel_POLK}.

We establish that the asymptotic bias of this method is a tunable constant depending on the compression parameter. These results yield an approach to importance reweighting that mitigates particle degeneracy, i.e., retaining a large number of particles with small weights \cite{li2014fight}, by directly compressing the embedding estimate of the posterior in the RKHS domain, rather than statistical tests that require sub-sampling in the distributional space \cite{bugallo2017adaptive,elvira2017improving}. The compression is performed \emph{online}, without waiting until the total number of samples $N$ are available, which is typically impractical.
Experiments demonstrate that this approach yields an effective tradeoff of consistency and memory, in contrast to the classical curse of dimensionality of MC methods. 

{\noindent \bf Additional Context.} Dimensionality reduction of nonparametric estimators been studied in disparate contexts. A number of works fix the sparsity dimension and seek the best $N$-term approximation in terms of estimation error. When a likelihood model is available, one terms the resulting active set a Bayesian \emph{coresets} \cite{campbell2019automated}. Related approaches called ``herding"  assume a fixed number of particles and characterize the resulting error in a kernel-smoothed density approximation \cite{chen2010super,lacoste2015sequential,keriven2018sketching,futami2019bayesian}. In these works, little guidance is provided on how to determine the number of points to retain.

 In contrast, dynamic memory methods automatically tune the number of particles to ensure small model bias. For instance, in \cite{campbell2018bayesian}, a rule for retainment based on gradient projection error (assuming the likelihood is available) is proposed, similar to those arising in kernel regression \cite{Koppel_POLK}. Most similar to our work is the setting where a likelihood/loss is unavailable and one must resort to metrics on density estimates, i.e., statistical tests, for whether new particles are significant. Specifically, in particle filters, multinomial resampling schemes can be used with Chi-squared tests to determine whether the current number of particles should increase/decrease \cite{pmlr-v32-jun14,elvira2016adapting}. The performance of these approaches have only characterized when their budget parameter goes to null or sparsity dimension goes to infinity. In contrast, in this work, we are especially focused on finite-sample analysis when budget parameters are left fixed, in order to elucidate the tradeoffs between memory and consistency, both in theory and practice.
 
%\vspace{-3mm}
%%%%%%%%%%%%%%%%%%%%%%%%%%%%%%%%%%%%%%%%%%%%%%%%%%%%%%%%%%%%%%%%%%%%%%%%%%%%%%%%%%%%%%%%%%%%%%%%%%%%%%%%%%%%%%%%%%%%%%%%%%%%%%%%%%%%%%%%%%%%%%%%%%%%%%%S  E  C  T  I  O  N %%%%%%%%%%%%%%%%%%%%%%%%%%%%%%%%%%%%%%%%%%%%%%%%%%%%%%%%%%%%%%%%%%%%%%%%%%%%%%%%%%%%%%%%%%%%%%%%%%%%%%%%
%%%%%%%%%%%%%%%%%%%%%%%%%%%%%%%%%%%%%%%%%%%%%%%%%%%%%%%%%%%%%%%%%%

\section{Elements of Importance Sampling}\label{sec:prob}

%\subsection{Importance Sampling}\label{subsec:importance_sampling}
%
In Bayesian inference \cite{sarkka2013bayesian}[Ch. 7], we are interested in computing expectations 
%\textcolor{blue}{In a lot of places we wrote $\bbx^n$ rather than $\bbx(n)$. Please double check to make sure I corrected them all.}
\begin{align}\label{eq:template_problem}
I(\phi):=\mathbb{E}_{  \bbx}[{\phi}(\bbx)\given \bby_{1:K}  ] = \int_{\bbx\in\ccalX} \phi(\bbx)  p(\bbx|\bby_{1:K})  d\bbx
\end{align}
on the basis of a set of available observations  {$ \bby_{1:K}:=\{\bby_{1:K}\}_{k=1}^{K}$}, where $\phi : \ccalX\rightarrow \reals$ is an arbitrary function, $\bbx$ is a random variable taking values in $\ccalX\subset \reals^p$ which is typically interpreted as a hidden parameter, and $\bby$ is some observation process whose realizations $\bby_{1:K}$ are assumed to be informative about parameter $\bbx$. {For example, $\phi(\bbx) = \bbx$ yields the computation of the posterior mean, and $\phi(\bbx)=\bbx^p$ denotes the $p$-th moment.}
In particular, define the posterior density\footnote{Throughout, densities are with respect to the Lebesgue measure on $\reals^p$.}
\begin{align}\label{eq:posterior}
p\left(\bbx \given \bby_{1:K}\right)  = \frac{p\left(\bby_{1:K} \given  \bbx  \right) p\left(\bbx\right)}{p\left( \bby_{1:K}\right)}.
\end{align}
We seek to infer the posterior \eqref{eq:posterior} with $K$ data points $\bby_{1:K}$ available at the outset. Even for this setting, estimating \eqref{eq:posterior} has unbounded complexity \cite{li2005curse,tokdar2010importance} {when the form of the posterior is unknown}. Thus, we prioritize efficient estimates of \eqref{eq:posterior} from an online stream of samples from an \emph{importance density} to be subsequently defined.
Begin by defining posterior $q(\bbx) $ and un-normalized posterior $\widetilde q(\bbx)$ as
%
%\begin{equation}\label{eq:pi_tilde}
% $$ q(\bbx) = {\widetilde q(\bbx) }/{Z},$$
% \end{equation}
 %
 %and 
 %
 \begin{align}
 	 \!\!\!\!q(\bbx) = {\widetilde q(\bbx) }/{Z} \; , \quad \!\widetilde q(\bbx){:=}\widetilde q(\bbx\given \bby_{1:K}){=}p\left(\bby_{1:K} \given  \bbx  \right) p\left(\bbx\right),
 \end{align} 
where $\widetilde q(\bbx)$ integrates to normalizing constant $Z {:=} p\left( \bby_{1:K}\right)$\footnote{Note that $q(\bbx)$ and $\widetilde q(\bbx)$ depend on the data $\{\bby_{1:K}\}_{k \leq K}$, although we drop the dependence to ease notation.}.
In most applications, we only have access to a collection of observations $\bby_{1:K}$ drawn from a static conditional density $p( \bby_{1:K} \given \bbx )$ and a prior for $p(\bbx)$. Therefore, the integral \eqref{eq:template_problem} cannot be evaluated, and hence one must resort to numerical integration such as Monte Carlo methods. 
%\red{We need to careful about the notation related to $\bby$ here. I would say that we drop the dependence on the data (we do nothing with the data, hence the sooner we drop $y_k$, the better).}
%
%{\begin{example}(Localization)
% fill in later \victor{we need this for ML/SP, probably, but not a statistics/probability journal}
%%\victor{It depends on the journal, but I am not sure we need this motivating example.}
%\end{example}}
%
In Monte Carlo, we approximate \eqref{eq:template_problem} by sampling. Hypothetically, we could draw $N$  samples $\bbx(n)\sim  q(\bbx)$ and estimate the expectation in \eqref{eq:template_problem} by the sample average
\begin{equation}\label{eq:template_sample_mean}
\mathbb{E}_{ q(\bbx)}[\phi(\bbx) ] \approx \frac{1}{N}\sum_{n=1}^N \phi(\bbx(n)),
\end{equation}
but typically it is difficult to obtain samples $\bbx(n)$ from posterior $ q(\bbx)$ of the {random variable}. To circumvent this issue, define the \emph{importance density}  $\pi(\bbx )$\footnote{In general, the importance density could be defined over any observation process  $\pi(\bbx \given \{\bby_{k}\} )$, not necessarily associated with time indices $k=1,\dots,K$. We define it this way for simplicity.} with the same (or larger) support as true density $q(\bbx)$, and multiply and divide by $\pi(\bbx )$ inside the integral \eqref{eq:template_problem}, yielding
\begin{equation}\label{eq:importance_distribution_factoring}
\int_{\bbx\in\ccalX} \phi(\bbx) q(\bbx) d\bbx = \int_{\bbx\in\ccalX}\frac{ \phi(\bbx)  q(\bbx)  } {\pi(\bbx)} \pi(\bbx )d\bbx,
\end{equation}
where the ratio ${{q}(\bbx)}/{\pi(\bbx)}$ is the Radon-Nikodym derivative, or unnormalized density, of the target ${q}$ with respect to the proposal $\pi$.  Then, rather than requiring samples from the true posterior, one may sample from the importance density $\bbx(n)\sim \pi(\bbx ),\;n=1,...,N$, and approximate \eqref{eq:template_problem} as %with the unnormalized IS (UIS) estimator:
\begin{align}\label{eq:UIS_estimator}
  I_N(\phi) := & \frac{1}{N}\sum_{n=1}^N \frac{ q(\bbx(n)) } {\pi(\bbx(n) )}\phi(\bbx(n)) \nonumber \\
 =&   \frac{1}{NZ}\sum_{n=1}^N g{(\bbx(n))}\phi(\bbx(n)), 
%\label{eq:UIS_estimator_2}
%
\end{align}
%
%where 
%
%\begin{equation}\label{eq:importance_weight}
%%
%g{(\bbx(n))}:= \frac{  q(\bbx(n)) }{\pi(\bbx(n))}
%%
%\end{equation}
where 
\begin{align}\label{eq:posterior1}
	g{(\bbx(n))}\equiv \frac{  \widetilde{q}(\bbx(n)) }{\pi(\bbx(n))},
\end{align} are the importance weights. 
%We note that in practice, we cannot calculate  $q(\bbx(n))$ since the target density $q(\x)$ is unknown and hence we calculate it using Bayes rule as follows:% likelihood $p\left(\{\bby_t\}_{t\leq T} \given  \bbx  \right) p\left(\bbx\right)$ and prior $p\left( \{\bby_t\}_{t\leq T}\right)$ as follows
%%
%\begin{equation}\label{eq:posterior1}
%%
% q(\bbx(n))   = \frac{ \widetilde q(\bbx(n))}{\int  \widetilde q(\bbx)  d\x}.
%%
%\end{equation}
%%
%%%
%%\begin{equation}\label{eq:posterior1}
%%	%
%%	q(\bbx(n))   = \frac{ \widetilde q(\bbx(n)) p\left(\{\bby_{1:K}\}_{k\leq K} \given  \bbx(n)  \right) p\left(\bbx(n)\right)}{\int  \widetilde q(\bbx)  p\left(\{\bby_{1:K}\}_{k\leq K} \given  \bbx  \right) p\left(\bbx\right) d\x}.
%%	%
%%\end{equation}
%%%
%Substituting \eqref{eq:posterior1} into \eqref{eq:importance_distribution_factoring}, we obtain 
%%
%\begin{equation}\label{eq:importance_weight2}
%%
%g{(\bbx(n))}\equiv \frac{\widetilde q(\bbx(n))}{\pi (\bbx(n))}.
%%
%\end{equation}
%%
Note that \eqref{eq:UIS_estimator} is unbiased, i.e., $\mathbb{E}_{\pi(\bbx)}[  I_N(\phi) ] = \mathbb{E}_{ q(\bbx)}[\phi(\bbx) ]$ {and consistent with $N$.} Moreover, its variance depends on the importance density $\pi(\bbx ) $ approximation of the posterior \cite{elvira2019generalized}.

%\victor{I would focus on uni-dimensional ${\phi}$ (i.e. without the boldface).}
% that's correct. The boldface is a typo -- on the previous page, g was defined as having real-valued scalar range
%[\red{Notation correction required, please follow the notation of (6-(8))}, oir I have used the correct notations in the Algorithms, for the theory, I was little confused]
%
Example priors and measurement models include Gaussian, Student's t, and uniform. Which choice is appropriate depends on the context \cite{sarkka2013bayesian}. % {another comment about this}.
%
%The challenge with IS is that we require evaluation of the posterior $ \mu(\bbx)$, but usually the normalizing constant $Z$ is not available, so in practice the evaluation is on the un-normalized density $\widetilde \mu(\bbx)$.
 The normalizing constant $Z$ can be also estimated with IS as 
\begin{equation}\label{eq:importance_sampling_denominator}
\hat Z := \frac{1}{N}\sum_{n=1}^N  g{(\bbx(n))}. 
\end{equation}
%
%%%%%%%%%%%%%%%%%%%%%%%%%%%%%%%%%%%%%%%%%%%%%%%%%%%%%%%%%%%%%%
%%%%%%%%%%%%%%%%%%%%%%%%%%%%%%%%%%%%%%%%%%%%%%%%%%%%%%%%%%%%%%
%%% A  L  G  O  R  I  T  H  M  %%%%%%%%%%%%%%%%%%%%%%%%%%%%%%%%%%%%%%%%%%%%%%%%%%%%%%%%%%%%%%
%%%%%%%%%%%%%%%%%%%%%%%%%%%%%%%%%%%%%%%%%%%%%%%%%%%%%%%%%%%%%%
%%%%%%%%%%%%%%%%%%%%%%%%%%%%%%%%%%%%%%%%%%%%%%%%%%%%%%%%%%%%%%\
\begin{algorithm}[t]
\caption{ IS: Importance Sampling with streaming samples}
\begin{algorithmic}
\label{alg:importance_sampling}
\REQUIRE Observation model  $p(\bby \given \bbx )$  and prior $p(\bbx)$ or target density $q(\x)$ (if known), importance density $\pi(\bbx)$. Set of observations $\{\bby_{1:K}\}_{k=1}^K$.% Define unnormalized version $\widetilde{q}(\bbx|\bby) \equiv \widetilde{q}(\bbx)$. 
\FOR{$N=0,1,2,\ldots$}
	%\STATE Observe $\bby_t$ to obtain latest observation model $p( \{\bby_u\}_{u\leq t} \given \bbx )$
	\STATE Simulate one sample from importance dist. $\bbx(n) \sim \pi(\bbx)$ \vspace{-0mm}
	\STATE Compute weight $g(\bbx(n))$ [cf. \eqref{eq:UIS_estimator}] 
	%if $q(\x)$ is known; else use \eqref{eq:importance_weight2}
%
%\begin{equation*}
%%
%g{(\bbx(n))}= \frac{ q(\bbx(n)) } {\pi(\bbx(n))} 
%%
%\end{equation*}
%
%\alec{Do we define $\mu(\bbx)$ anywhere in the algorithm pseudo-code? By $\mu(\bbx)$, we really mean $\mu(\bbx \given \bby)$ right?}
%\victor{Good question. I would drop the dependence of $\bbx$, and always denote it as $\mu(\bbx)$: note that $\mu(\bbx)$ is not a distribution since it does not integrate up to one, but $\mu(\bbx|\bby)$ gives the impression of conditional distribution. Since the data are always fixed, I would drop $\bby$, but I let you decide.}
%\alec{A few notational details to prioritize clarity:
%%
%\begin{enumerate}
%%
%\item I think we should have it be $\mu(\bbx|\bby)$ always to emphasize dependence on observations. Or if not, we need to specifically define the short-hand notation.
%%
%\item It seems more natural to let $\mu(\bbx|\bby)$ be the exact posterior and $\widetilde{\mu}(\bbx\given\bby)$ or $\hat{\mu}(\bbx\given\bby)$ be the unnormalized version. Intuitively to me the true distribution shouldn't have any adornment.
%%
%\item The importance distribution $\pi(\bbx)$ is really $\pi(\bbx\given \bby)$ right? Again, we need to define the short-hand notation.
%%
%\end{enumerate}
%%
%}
	\STATE Compute normalized weights $\overline{w}(n)$ by dividing by estimate for summand \eqref{eq:importance_sampling_denominator}: \vspace{-2mm}
\begin{equation*}
\overline{w}(n)= \frac{g{(\bbx(n))} }{\sum_{u=1}^N g{(\bbx(u))}}\ \ \ \ \ \text{for \ \ all}\ \ \ n. 
\end{equation*}\vspace{-2mm}
%\victor{I think we need to re-normalized all available weights, i.e., state that it is for $n=0,...,n$ (not only the $n$th)}
%
%
	\STATE Estimate the expectation with the self-normalized IS as  \vspace{-2mm}%conditional distribution
\begin{equation*}
I_N(\phi)= \sum_{n=1}^N \overline{w}(n) \phi(\bbx(n))
\end{equation*}\vspace{-2mm}
{	\STATE The posterior density estimate is given by\vspace{-2mm}
$$ {\mu}_N = \sum_{n=1}^N \overline{w}(n) \delta_{\bbx(n)}$$}\vspace{-2mm}
\ENDFOR
%\ENDALOOP
%\victor{Great!}
\end{algorithmic}
\end{algorithm}
%The approximation \eqref{eq:UIS_estimator} is always unbiased,  %likelihood $p(\bbx)$.
%
Note that in Eq. \eqref{eq:UIS_estimator}, the unknown $Z$ can be replaced by $\hat Z$ in \ref{eq:importance_sampling_denominator}. Then, the new estimator is given by
\begin{align}\label{eq:SNIS_estimator}
 I_N(\phi) :=&  \frac{1}{N\hat Z}\sum_{n=1}^N g{(\bbx(n))}\phi(\bbx(n)) \\
=&  \frac{1}{\sum_{j=1}^N g{(\bbx(j))}} \!\!\sum_{n=1}^N g{(\bbx(n))}\phi(\bbx(n))\nonumber
\\
=&\sum_{n=1}^N \overline w(n) \phi(\bbx(n)),\nonumber
\end{align}
where the { normalized }$\overline{w}(n)$ weights are defined 
\begin{equation}\label{eq:normalized_importance_weight}
\overline{w}(n)\equiv\frac{g{(\bbx(n))} }{\sum_{u=1}^N g{(\bbx(u))}},
\end{equation}
for all $n$. The whole IS procedure is summarized in Algorithm \ref{alg:importance_sampling}. The function $I_N(\phi)$ is the normalized importance sampling (NIS) estimator. It is important to note that the estimator $I_N(\phi)$ can be viewed as integrating a function $\phi$ with respect to density $ {\mu}_N$ defined as %\vspace{-2mm}
\begin{align}\label{eq:importance_sampling_empirical_measure}
 {\mu}_N(\x):=\sum\limits_{n=1}^{N}\overline{w}(n)\delta_{\bbx(n)},
\end{align}
which is called the particle approximation of $ {q}$. 
Here $\delta_{\bbx(n)}$ denotes the Dirac delta measure evaluated at  $\bbx(n)$. This delta expansion is one reason importance sampling is also referred to as a histogram filter, as they quantify weighted counts of samples across the space. Subsequently, we leave the argument (an event, or measurable subset) of the delta $\delta_{\bbx(n)}$ as implicit.

As stated in \cite{first,second,agapiou2017importance}, for consistent estimates of \eqref{eq:template_problem}, we require that $N$, the number of samples $\bbx(n)$ generated from the importance density, and hence the parameterization of the importance density, to go to infinity $N\rightarrow \infty$. Therefore, when we generate an infinite stream of particles, the parameterization of the importance density grows unbounded as it accumulates every particle previously generated.  {We are interested in allowing $N$, the number of particles, to become large (possibly infinite), while the importance density's complexity is moderate, thus overcoming an instance of the curse of dimensionality in Monte Carlo methods.} 
%
%
%The question, then, is how to obtain \emph{nearly} consistent posterior representations for \eqref{eq:template_problem} of finite memory. 
%
In the next section, we propose a method to do so.

\section{Compressing the Importance Distribution}\label{sec:algorithm}
%
%\red{Present the measure update with kernel and mention that the number of particles go to infinity}
%
%\red{introduce kernel mean embedding}
%
%\red{replace the original dictionary with the projected one and present updates}
%
In this section, we detail our proposed sequential compression scheme for overcoming the curse of dimensionality in importance sampling. However, to develop such a compression scheme, we first rewrite importance sampling estimates in vector notation to illuminate the dependence on the number of past particles generated. Then, because directly defining projections over measure spaces is intractable, we incorporate a distributional approximation called a mean embedding \cite{sriperumbudur2010hilbert}, over which metrics can be easily evaluated. This permits us to develop our main projection operator.

 Begin by noting the curse of dimensionality in importance sampling can be succinctly encapsulated by rewriting the last step of Algorithm \ref{alg:importance_sampling} in vector notation. Specifically, define  ${\bbg}_n \in \reals^n$  ,
%
%My proposal is, after each new sample $\bbx^{n}$, we have a (kernel) density estimate as in \eqref{eq:importance_sampling_empirical_measure} where the empirical summand goes from $n'=1,\dots,n$ just as it does in the last step of Algorithm \ref{alg:importance_sampling}. We can write Algorithm \ref{alg:importance_sampling} in a vector format with
%
 $\overline{\bbw}_n \in \reals^n$ and $\bbX_n = [\bbx(1) ; \cdots \bbx(n)] \in \reals^{p \times n}$. Then, after each new sample $\bbx(n)$ is generated from the importance distribution, we incorporate it into the empirical measure \eqref{eq:importance_sampling_empirical_measure} through the parameter updates 
%
% \cred{
% \begin{equation}\label{eq:importance_sampling_stack}
% %
% {w}(n)= \frac{ \mu(\bbx^{n}) } {\pi(\bbx^{n})}  \; , \quad \bbw_{n+1} = [\bbw_n ; g{(\bbx^{n})} ] \;, \quad \overline{\bbw}_{n+1} = \frac{1}{\sum_{n'=1}^n w^{(n')}}\bbw_{n+1} \;, \qquad \bbX_{n+1} = [\bbX_n ; \bbx^{n} ] 
% %
% \end{equation}
% }
%\victor{I have slightly corrected the notation. I think it will be more practical as it is now:}
\begin{align}\label{eq:importance_sampling_stack}
\bbg_{n} =& [\bbg_{n-1} ; g{(\bbx(n))} ] \;,\nonumber  \\
 \quad \overline{\bbw}_{n} =& z_n \bbg_{n} \;, \qquad \bbX_{n} = [\bbX_{n-1} ; \bbx(n) ],
\end{align}
 where we define $z_n:=1/(\mathbf{1}_n^T\bbg_{n})$ and  $\mathbf{1}_n$ is the all ones column vector  with dimension $n$.  The unnormalized posterior density estimate parameterized by $\bbg_n$ and \emph{dictionary} $\bbX_n$ is given by
\begin{align}\label{eq:posterior_distribution}
 \tilde{\mu}_n = \sum_{u=1}^n \bbg_n(u) \delta_{\bbx(u)},
\end{align}
%\victor{I would define here $z_n:=1/(\sum_{n'=1}^n w^{(n')})$, since above it is only described for $N$.}
where we define $\bbg_n(u) := g(\bbx(u))$ is the importance weight \eqref{eq:UIS_estimator} and $\delta_{\bbx(u)}$ is the Dirac delta function, both evaluated at sample $\bbx(u)$. Denote as $\Omega_{\bbX_n}$ the measure space ``spanned" by Dirac measures centered at the samples stored in dictionary $\bbX_n$. More specifically, given measurable space $(\Omega,\Sigma,\lambda)$, where $\Omega$ is a set of outcomes, $\Sigma$ is a $\sigma$-algebra whose elements are subsets, and $\lambda: \Sigma \rightarrow \reals$ denotes the Lebesgue measure, we define the restricted $\sigma$-algebra as $\Sigma_{\bbX_n}=\{F\cap\{\bbx(u)\}_{u\leq n} : F\in\Sigma\}$. The measure space associated with $\Sigma_{\bbX_n}$ and the Lebesgue measure $\lambda_{\bbX_n}$ over this restricted $\sigma$-algebra we denote in shorthand as $\Omega_{\bbX_n}$ (see \cite{wheeden1977measure} for more details).

We note that the unnormalized posterior density in \eqref{eq:posterior_distribution} is a linear combination of (nonnegative) Dirac measures with mass $\bbg_n(u)$ for each sample $\x(u)$. 
%Hence, $\tilde{\mu}_n\in\Omega_{\bbX_n}$ defines a signed measure  where $\bbg_n(u)$ is allowed to take real values.
% \green{Victor: I am not sure if I see why the measure is signed in the alg. Then, the 'g' can be negative? Does it happen in reality? Do we ensure that the (unnormalized or normalized) kernel approximation of the density is non-negative in all the domain (i.e., for any subset)? \blue{Amrit: Thanks for the comments. This was a mistake earlier. $\mu$ here is not a signed measure because $\bbg_n$ is always positive (it is the rations of probability distributions.) I have updated the discussion and removed the associated statement.}}
 The question is how to select a subset of columns of $\bbX_n$ and {modify} the weights $\bbg_n$ such that with an infinite stream of $\bbx(n)$, we  ensures the number of columns of $\bbX_n$ is finite and the empirical integral estimate tends to its population counterpart, i.e., the integral estimation error becomes very small (or goes to zero) as $n$ tends to infinity \cite{ghosal2000convergence}. Henceforth, we refer to the number of columns in matrix  $\bbX_n$ which parameterizes \eqref{eq:posterior_distribution} as the \emph{model order} {denoted by $M_n$}.
%\red{I talked about the definition of measure space in the previous paragraph because we need to mention measure space before introducing the kernel mean embedding which is defined with respect to a measure space. We might want to talk about it by stating things more rigorously like  and $\Omega_{\bbX_n}$ defines the set of Dirac measures on this space.}	

%We will use the notation $\tilde{\mu}_n (\bbg_n,\bbX_n) $ to emphasize the fact that the estimate $\tilde{\mu}_n$ is parameterized by the weights $\bbg_n$ and the dictionary $\bbX_n$. 
%
\begin{figure}[t] 
	\centering 
	\includegraphics[scale=0.25]{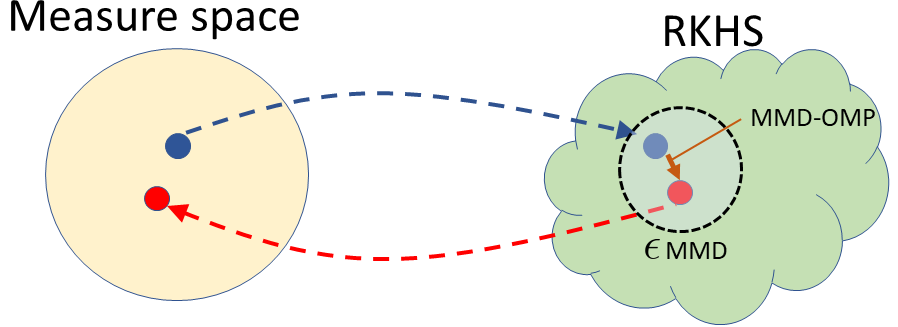}
	\caption{Approximating the distributions via kernel mean embedding.} 
	\label{diag}\vspace{-6mm}
\end{figure}
\subsection{Kernel Mean Embedding}
In this subsection, we introduce kernel mean embedding which maps the measure estimate in the measure space $\Omega_{\bbX_n}$ to the corresponding reproducing kernel Hilbert space (RKHS) denoted by $\mathcal{H}_{{\bbX_n}}$ as shown in Fig. \ref{diag}. 	There is a kernel associated with with RKHS defined as $\kappa:\mathcal{X}\times\mathcal{X}\rightarrow\mathbb{R}$ and $\kappa(\bbx,\cdot)\in\mathcal{H}_{\bbX_n}$.  In order to appropriately select a subset of $\bbX_n$ which would define an approximator for $\tilde{\mu}_n$ in \eqref{eq:posterior_distribution},  we employ the notion of kernel mean embedding \cite{sriperumbudur2010hilbert} of Dirac measures, and then perform approximation in the corresponding RKHS $\mathcal{H}_{\bbX_n}$. 
Doing so is motivated by the fact that operations involving distributions embedded in the RKHS may be evaluated in closed form, whereas in general measure spaces it is often intractable. The explicit value of the map for $ \tilde{\mu}_n$ to RKHS is given by 
	 \begin{align}\label{eq:kernel_mean_emb}
		\beta_{n}=\sum_{u=1}^{n} \bbg_n(u) \kappa(\bbx(u),\cdot),
	\end{align}
	 and $\beta_{n}\in\mathcal{H}_{\bbX_n}$. We remark here that the associated kernel $\kappa$ with RKHS is assumed to be a  characteristic kernel which ensure that the mapping $\tilde{\mu}_n\rightarrow\beta_{n}$ is injective \cite[Def. 3.2]{kernel_mean_embedding}. A characteristic  kernel is imperative to make sure that $\|\beta_{\mathbb{P}}-\beta_{\mathbb{Q}}\|_{\mathcal{H}}=0$ if and only if ${\mathbb{P}}={\mathbb{Q}}$ for measures ${\mathbb{P}}$ and ${\mathbb{Q}}$, i.e., that it satisfies the identity of indiscernibles. This makes sure that there is no information loss by introducing the mapping via kernel mean embedding.  
	 
For practical purpose, we are interested in obtaining the value of the underlying posterior density  associated with the mean embedding, which necessitates a way to invert the embedding. Doing so is achievable by computing the distributional pre-image \cite{honeine2011preimage}, which for a given $\beta_{n}\in\mathcal{H}_{{\bbX}_n}$ is given by 
	 \begin{align}\label{eq:proximal_hilbert_representer}
		\bbg^\star : =&\argmin_{\bbg \in \reals^{n}} \!
		\Big\| \beta_{{n}}
		\!\! -\! \!\sum_{u=1}^{n}{\bbg}(u)\kappa{(\bbx(u),\cdot)} \Big\|_\mathcal{H}^2,
	\end{align}
%\green{Victor: by looking at (14), is not $g^*$ exactly $g_n$? Maybe I am missing something (or maybe the number of elements should not be the same in both cases?)} [\blue{Amrit: Yeah, you are right, I tried to explain it below in blue. }]
where $\sum_{u=1}^{n}{\bbg}(u)\kappa{(\bbx(u),\cdot)} $ is the kernel mean embedding for the Dirac measure $\sum_{u=1}^{n}{\bbg}(u)\delta_{\bbd_{\bbx(u)}}$ and ${\bbg}=\bbg^\star$ is obtained as the solution of \eqref{eq:proximal_hilbert_representer}. Therefore, for a given  $\beta_{n}\in\mathcal{H}_{\bbX_n}$, we could recover the corresponding distribution in the measure space $\Omega_{\bbX_n}$  by solving \eqref{eq:proximal_hilbert_representer}. Note that \eqref{eq:proximal_hilbert_representer} exhibits a closed form solution with $\bbg^\star=\bbg_n$, as the distributional measures have a Dirac measure structure. This motivates us to perform the compression in the RKHS $\mathcal{H}_{\bbX_n}$ parameterized by $\bbX_n$.    
% %, we could write }
%	%
%	 \blue{\begin{align}\label{eq:proximal_hilbert_representer2}
%		\bbg^\star : =&\argmin_{\bbg \in \reals^{|{\bbX}_n|}} \!
%		\Big\| \sum_{u=1}^{|{\bbX}_n|} \bbg_n(u) \kappa(\bbx(u),\cdot)
%		\!\! -\! \!\sum_{u=1}^{|{\bbX}_n|}{g}(u)\kappa{(\bbd_u,\cdot)} \Big\|_\mathcal{H}^2
%		%
%	\end{align}}
%
 		
	Specifically, given past particles collected in a dictionary $\bbX_n$, we seek to select the subset of columns of $\bbX_n$ to formulate its compressed variant $\bbD_n$. We propose to project the kernel mean embedding $\beta_{n}$ onto subspaces $\mathcal{H}_{\bbD_n}=\text{span}\{\kappa(\bbd_u,\cdot) \}_{u=1}^{M_n}$ that consist only of functions that can be represented using dictionary $\bbD_n = [\bbd_1,\ \ldots,\ \bbd_{M_n}] \in \reals^{p \times {M_n}}$. 
		More precisely, $\ccalH_{\bbD_n}$ is defined as a subspace in the RKHS $\mathcal{H}_{\bbX_n}$ that can be expressed as a linear combination of kernel evaluations at points $\{\bbd_u\}_{u=1}^{M_n}$.  We enforce efficiency by selecting dictionaries $\bbD_n$ such that $M_n \ll n$. The eventual goal is to achieve a finite memory $M_n=\mathcal{O}(1)$ as $n\rightarrow\infty$ or $\frac{M_n}{n}\rightarrow 0$ as $n\rightarrow\infty$, while ensuring the underlying integral estimate has minimal bias, i.e., that we can obtain nearly consistent finite particle estimates. 

%%%%%%%%%%%% %%%%%%%%%%%%%%%%%%%%%%%%%%%%%%%%%%%%%%%%%%%%%%%%%%%%%%%%%%%%%%
%%%%%%%%%%%%%%%%%%%%%%%%%%%%%%%%%%%%%%%%%%%%%%%%%%%%%%%%%%%%%%%%%%%%%%%%%%
%%%%%%%%%%%%%%% A  L  G  O  R  I  T  H  M  %%%%%%%%%%%%%%%%%%%%%%%%%%%%%%%%%%%%%%%%%%%%%%%%%%%%%%%%%%%%%%
%%%%%%%%%%%%%%%%%%%%%%%%%%%%%%%%%%%%%%%%%%%%%%%%%%%%%%%%%%%%%%%%%%%%%%%%%%
%%%%%%%%%%%%%%%%%%%%%%%%%%%%%%%%%%%%%%%%%%%%%%%%%%%%%%%%%%%%%%%%%%%%%%%%%%\
\begin{algorithm}[t]
\caption{Compressed Kernelized IS (CKIS)}
\begin{algorithmic}
\label{alg:compressed_online_kernel_importance_sampling}
\REQUIRE Unnormalized target distribution $\widetilde q(\x)$, importance distribution $\pi(\bbx)$.%, Observation collection $\{\bby_k\}_{k=1}^K$%, and initial unnormalized posterior $ \tilde{\mu}_n$ at $n=0$..
\FOR{$n=0,1,2,\ldots,N$}
	%\STATE Observe $\bby_n$ to obtain latest observation model $p( \{\bby_u\}_{u\leq t} \given \bbx )$
	\STATE Simulate one sample from importance dist. $\bbx(n) \sim \pi(\bbx)$
	\STATE Compute the importance weight $g{(\bbx(n))}\equiv \frac{  \tilde q(\bbx(n)) }{\pi(\bbx(n))}$
%\begin{equation*}
%%
%g(\bbx^{n})= \frac{ q(\bbx^{n}) } {\pi(\bbx^{n})} 
%%
%\end{equation*}
%
	\STATE Normalize weights $w(n)$ by estimate for summand \eqref{eq:importance_sampling_denominator}: \vspace{-2mm}
\begin{equation*}
\overline w(j)\!:=\! \frac{{w}(j) }{z_n},  j=1,...,n. \; ,  z_n=\sum_{u=1}^n {w}(u)
\end{equation*}\vspace{-3mm}
	\STATE Update the mean embedding via last sample $\&$ weight [cf. \eqref{eq:kernel_mean_emb2}] \vspace{-2mm}
	$$	\tilde\beta_{n}=\beta_{{n-1}} + \bbg_n(n) \kappa(\bbx(n),\cdot). $$
	
	\STATE Append dictionary $\tbD_{n} \!= \![\bbD_{n-1} ; \bbx(n)] $ and importance weights $\tilde \bbg_{n} \!= \! [\bbg_{n-1} ; g(\bbx(n)) ] $ 
	\STATE Compress the mean embedding as (cf. Algorithm \ref{alg:komp}) \vspace{-2mm}
	$$(\beta_{n},\bbD_{n},\bbg_{n}) = \textbf{MMD-OMP}(\tilde\beta_{n},\tbD_{n},\tbg_{n},\epsilon_n)$$

	%\STATE \alec{this last step could be omitted/explained as the output of running the algorithm}
  \STATE Evaluate the pre-image to calculate $\hat\mu_n$ using \eqref{eq:proximal_hilbert_representer}

\STATE  Estimate the expectation as \vspace{-0mm}%conditional distribution $p(\bbx \given \{\bby_u\}_{u\leq t}) $ as \eqref{eq:importance_sampling_estimator}
%\begin{equation*}
%
$\hat{I}_n= \sum_{u=1}^{|\bbD_{n}|}\overline  w(u) \phi(\bbx(u))$ %\textcolor{blue}{?????????? \quad  \text{please correct/update}}
%
%\end{equation*}\vspace{-4mm}
%
\ENDFOR
%\ENDALOOP
\end{algorithmic}
\end{algorithm}

\subsection{Greedy Subspace Projections}
Next, we develop a way to only retain statistically significant particles by appealing to ideas from subspace based projections \cite{rockafellar1976monotone}, inspired by \cite{Koppel_POLK,koppelconsistent}. 
 To do so, we begin by rewriting the evolution of the mean embedding in \eqref{eq:kernel_mean_emb} as
%
%, its weight and dictionary, are instead updated as \eqref{eq:importance_sampling_stack}, which may be rewritten as
%
%
%
 	\begin{align}\label{eq:kernel_mean_emb2}
	\beta_{n}=\beta_{{n-1}} + \bbg_n(n) \kappa(\bbx(n),\cdot).
\end{align}
 We note that the update in \eqref{eq:kernel_mean_emb2} can be written as 
%	%
	\begin{align}\label{eq:proximal_hilbert_dictionary}
\beta_{n}=& 
%\argmin_{f\in\mathcal{H}}
%\ \|f-\big(\beta_{{\mu}_{n-1}} + \bbg_n(n) \kappa(\bbx(n),\cdot)\big)\|^2_\mathcal{H} \nonumber\\
%=&
 \argmin_{f\in\mathcal{H}_{\bbX_n}} \ \|f-\big(\beta_{{n-1}} + \bbg_n(n) \kappa(\bbx(n),\cdot)\big)\|^2_\mathcal{H} \; , 
\end{align}
where the  equality employs the fact that $\beta_{n}$ can be represented using only the elements in $\mathcal{H}_{\bbX_n}=\text{span}\{\kappa(\bbx(u),\cdot)\}_{u\leq n}$. 
 Observe that \eqref{eq:proximal_hilbert_dictionary} defines a projection of the update $\big(\beta_{\tilde{\mu}_{n-1}} + \bbg_n(n) \kappa(\bbx(n),\cdot)\big)$ onto the subspace defined by $\mathcal{H}_{\bbX_n}$, which we propose to replace at each iteration by a projection onto a subspace defined by dictionary $\bbD_{n}$, which is extracted from the particles observed thus far. The process by which we select $\bbD_{n}$ will be discussed next. To be precise, we replace the update \eqref{eq:proximal_hilbert_dictionary} in which the number of particles grows at each step by the subspace projection onto $\mathcal{H}_{\bbD_{n}}$ as 
\begin{align}\label{eq:proximal_hilbert_dictionary2}
	\beta_{n}=& \argmin_{f\in\mathcal{H}_{\bbD_n}}
	\ \|f-\big(\beta_{\tilde{\mu}_{n-1}} + \bbg_n(n) \kappa(\bbx(n),\cdot)\big)\|^2_\mathcal{H} \nonumber
	\\
	:=& \mathcal{P}_{\mathcal{H}_{\bbD_n}}[\beta_{{n-1}} + \bbg_n(n) \kappa(\bbx(n),\cdot)].
%	=& \mathcal{P}_{\mathcal{H}_{\bbD_n}}[\tilde\beta_{{\mu}_n}] \; , 
	%
	\end{align}
Let us define $\tilde\beta_{n}:=\beta_{{n-1}} + \bbg_n(n) \kappa(\bbx(n),\cdot)$, which means that $\beta_{n}= \mathcal{P}_{\mathcal{H}_{\bbD_n}}[\tilde\beta_{n}]$.
	Let us denote the corresponding dictionary update as 
	\begin{align}\label{eq:importance_sampling_stack2}
		\tbg_{n} =& [\bbg_{n-1} ; g{(\bbx(n))} ] \;, \qquad \tbD_{n} = [\bbD_{n-1} ; \bbx(n) ] \nonumber  \\
		\quad \overline{\bbw}_{n} =& z_n \bbg_{n} \;, 
	\end{align}
where $\bbD_{n-1}$ has $M_{n-1}$ number of elements and $\tbD_{n}$ has $\tilde M_n= M_{n-1}+1$. Using the expression for mean embedding in \eqref{eq:kernel_mean_emb}, we may write the projection in \eqref{eq:proximal_hilbert_dictionary2} as follows
	\begin{align}\label{eq:proximal_hilbert_representer22}
	\bbg_{n}:=&\argmin_{\bbg \in \reals^{{M}_{n}}}
	\Big\|\! \sum_{s=1}^{{M}_{n}} {\bbg}(s)\kappa{(\bbd_s,\cdot)}
	\! - \!\sum_{u=1}^{\tilde{M}_n}\tilde{\bbg}_n(u)\kappa{(\tbd_u,\cdot)} \Big\|^2_\mathcal{H}\\
	%
	%
	%\! \argmin_{\bbg \in \reals^{{M}_{n}}} \! 
	% \! \Big( \! \sum_{s,{s'}=1}^{{M}_{n}} \!\!\!\! {g}_s{g}_{s'} \delta_{\bbd_{s}}\delta_{\bbd_{s}'}
	%  \!-\! 2 \! \! \sum_{s,u=1}^{{M}_{n},\tilde{M}_n} \! \!\!\!{g}_s\tilde{g}_u \delta_{\bbd_s}\delta_{\tbd_u} \!\!\!+\!  \!\!\!\!  \sum_{u,{u'}=1}^{\tilde{M_n}} \!\!\!\!\tilde{g}_{u}\tilde{g}_{u'}\delta_{\tbd_{u}}\delta_{\tbd_{u'}}\! \!  \Big)  \nonumber \\
	%
	&\hspace{-8mm}=\argmin_{\bbg \in \reals^{{M}_{n}}} \! 
	\left( \!\bbg^T\bbK_{\bbD_{n},\bbD_{n}}\bbg
	\!-\! 2\bbg^T\!\bbK_{\bbD_{n},\tbD_{n}}\tbg_{n}
	\!+\! \tbg_{n}\bbK_{\tbD_{n},\tbD_{n}}\tbg_{n} \!\right) \!
	\!\nonumber \; ,
\end{align}
	where we expand the square and define the kernel covariance matrix $\bbK_{\bbD_{n},\tbD_{n}}$ whose $(s,u)^\text{th}$ entry is given by $\kappa{(\bbd_s,{\tbd_u})}$. The other matrices $\bbK_{\tbD_{n},\tbD_{n}}$ and $\bbK_{\bbD_{n},\bbD_{n}}$ are similarly defined. The problem in \eqref{eq:proximal_hilbert_representer22} may be solved explicitly by computing gradients and solving for $\bbg_{n}$ to obtain
\begin{equation} \label{eq:hatparam_update}
	\bbg_{n}=  \bbK_{\bbD_{n},\bbD_{n}}^{-1}\bbK_{\bbD_{n} \tbD_{n}} \tbg_{n}. \; 
\end{equation}
Then, the projected estimate of the mean embedding $\tilde\beta_{n}$ is given by 
 	\begin{align}
		\beta_{n}=& \mathcal{P}_{\mathcal{H}_{\bbD_n}}[\tilde\beta_{n}]
		=\sum_{s=1}^{{M}_{n}} {\bbg_n}(s)\kappa{(\bbd_s,\cdot)},
 	\end{align}
	where $\bbg_n$ is obtained as a solution to \eqref{eq:hatparam_update}. Now, for a given dictionary $\bbD_n$, we know how to obtain the projected version of the mean embedding for each $n$. Next, we discuss the procedure to obtain $\bbD_n$ at each $n$.

\begin{algorithm}[t]
	\caption{MMD based Orthogonal Matching Pursuit (MMD-OMP) \hspace{-2mm}}
	\label{alg:komp}
	\begin{algorithmic}
		\REQUIRE  kernel mean embedding $\tilde\beta_{n}$ defined by dict. {$\tbD_n \in \reals^{p \times \tilde{M}_n}$}, coeffs. {$\tbg \in \reals^{\tilde{M}_n}$}, approx. budget  $\epsilon_n > 0$ \\
		\textbf{initialize} $\beta=\tilde\beta_{n}$, dictionary $\bbD = \tbD_n$ with indices $\ccalI$, model order $M=\tilde{M}_n$, coeffs.  $\bbg = \tbg_n$.
		\WHILE{candidate dictionary is non-empty $\ccalI \neq \emptyset$}
		\FOR {{$j=1,\dots,{M}$}}
			\STATE Find minimal approx. error with dictionary element $\bbd_j$ removed \vspace{-2mm}
			{			$$\gamma_j =  \text{MMD}\Big[\tilde\beta_{n}, \sum_{k \in \ccalI \setminus \{j\}} \bbg(k) {\kappa}(\bbd_k, \cdot)\Big] \; .$$ \vspace{-5mm}}
			\ENDFOR
		\STATE Find dictionary index minimizing approximation error: $j^\star = \arg\min_{j \in \ccalI} \gamma_j$
		%    \State Add dictionary element $\tbd_{\gamma}$ to basis $\bbD \leftarrow [\bbD,\;\;\tbd_{\gamma}]$; increment model order ${M} \leftarrow {M}+1$.
		\INDSTATE{{\bf{if }}  minimal approx. error exceeds threshold $\!\gamma_{j^\star}\!\! > \! \epsilon_n$}
		\INDSTATED{\bf stop} 
		\INDSTATE{\bf else} 
		
		\INDSTATED Prune dictionary $\bbD\leftarrow\bbD_{\ccalI \setminus \{j^\star\}}$, {remove the columns associated with index $j^\star$}
		\INDSTATED Revise set $\ccalI \leftarrow \ccalI \setminus \{j^\star\}$ and {model order ${M} \leftarrow {M}-1$.}
		\INDSTATED Update weights $\bbg$ defined by current dictionary $\bbD$
		%    %
		\vspace{-2mm}{$$\bbg = {\argmin_{\bbw \in \reals^{{M}}}} \big\| \tilde\beta_{n} - \sum_{u=1}^{M} \bbw(u) {\kappa}(\bbd_u, \cdot) \big\|_{\ccalH}$$}\vspace{-3mm}
		%
	%	\red{$M_n$ (without the $\tilde$) is not present in this algorithm. Please unify the notation with the discussion of previous page and the dictionary update subsection below (or maybe I'm missing a detail here! :-) )}
		\INDSTATE {\bf end}
		\ENDWHILE

		Assign $\bbg_n=\bbg$ and {Evaluate} the projected kernel mean embedding as 
		%	\begin{align}
			$\beta_{n}=\sum_{u=1}^{{M}} {\bbg_n}(u)\kappa{(\bbd_u,\cdot)}$
		%\end{align}}
		\RETURN $\beta_{n},\bbD_n,\bbg_n$ of complexity $\!M\!\! \leq  \!\tilde{M}$ s.t. $\text{MMD}[\tilde\beta_{n},\beta_{n}]\leq \!\eps_n$ 
	\end{algorithmic}
\end{algorithm}

%%%%%%%%%%%%%%%%%%%%%%%%%%%%%%%%%%%%%%%%%%%%%%%%%%%%%%%%%%%%%%
%%   S  U  B  S  E  C  T  I  O  N    %%%%%%%%%%%%%%%%%
%%%%%%%%%%%%%%%%%%%%%%%%%%%%%%%%%%%%%%%%%%%%%%%%%%%%%%%%%%%%%%
\subsection{Dictionary Update}  The selection procedure for the dictionary $\bbD_{n}$ is based upon greedy sparse approximation, a topic studied in compressive sensing  \cite{needell2008greedy}.  The function $\tilde\beta_{n}:=\beta_{{n-1}} + \bbg_n(n) \kappa(\bbx(n),\cdot)$ is parameterized by dictionary $\tbD_{n}$, whose model order is $\tilde{M}_n={M}_{n-1} +1$. We form $\bbD_{n}$ by selecting a subset of $M_{n}$ columns from $\tbD_{n}$ that are best for approximating the kernel mean embedding $\tilde\beta_{_n}$ in terms of maximum mean discrepancy (MMD). As previously noted, numerous approaches are possible for sparse representation. We make use of destructive \emph{ orthogonal matching pursuit} (OMP) \cite{Vincent2002} with allowed error tolerance $\epsilon_n$ to find a dictionary matrix $\bbD_{n}$ based on the one that includes the latest sample point $\tbD_{n}$. With this choice, we can tune the stopping criterion to guarantee the per-step estimates of mean embedding are close to each other. We name the compression procedure {MMD-OMP} and it is summarized in Algorithm \ref{alg:komp}.
From the procedure in Algorithm \ref{alg:komp}, note that the projection operation in \eqref{eq:proximal_hilbert_dictionary2} is performed in a manner that ensures that   $\text{MMD}[\tilde\beta_{n},\beta_{n}]\leq \eps_n$ for all $n$, and we recall that $\beta_{n}$ is the compressed version of $\tilde\beta_{n}$.

\section{Balancing Consistency and Memory}\label{sec:convergence}
%
%\textcolor{blue}{In a lot of places we wrote $\bbx^n$ rather than $\bbx^{(n)}$. Please double check to make sure I corrected them all, especially in the proofs in appendices. \\  \\
%
%\textcolor{blue}{Also please check that I didn't screw up the difference between $\mu$, $\pi$, $q$, and $m$ and their respective tilde/hat versions.
%}

In this section, we characterize the convergence behavior of our posterior compression scheme. Specifically, we establish conditions under which the asymptotic bias is proportional to the kernel bandwidth and the compression parameter using posterior distributions given by Algorithm \ref{alg:compressed_online_kernel_importance_sampling}. To frame the discussion, we note that the NIS estimator \eqref{eq:SNIS_estimator} $I_N(\phi)$, whose particle complexity goes to infinity, is asymptotically consistent \cite{mcbook}[Ch. 9, Theorem 9.2], %, i.e.,
%
%%
%\begin{align}\label{eq:consistent}
%%
% \mathbb{P}\left(\lim_{N\rightarrow \infty}I_N(\phi)= I(\phi) \right)=1 .%, \quad \Var(\hat{I})=\mathbb{E}\left[  \right]
%\end{align}
%%
and that the empirical posterior $\mu_N(\cdot)$ contracts to its population analogue at a $\ccalO(1/N)$ rate where $N$ is the number of particles. To establish consistency, we first detail the technical conditions required. 
%detailed in Appendix \ref{apx_technicalities} are required.
\subsection{Assumptions and Technical Conditions}\label{apx_technicalities}
%%%%%%%%%%%%%%%%%%%%%%%%%%%%%%%%%%%%%%%%%%%%%%%%%%%%%%%%%%%%%%%%%%%%%%%%%%%%%%%%%%%%%%%%%%%%%%%%%%%%%%%%%%%%%%%%%%%%%%%%%%%%%%%%%%%%%%%%%%%
\begin{assumption}\label{as:first}
Recall the definition of the target distribution $q$ from Sec. \ref{sec:prob} (following \eqref{eq:posterior}). Denote the integral of test function $\phi:\ccalX \rightarrow \reals$ as $q(\phi)$. 
\begin{enumerate}[label=(\roman*)]
\item Assume that $\phi$  is absolutely integrable, i.e., $q(|\phi|)<\infty$, and has absolute value at most unit $|\phi|\leq 1$. \label{as:testfunction1}
\item The test function has absolutely continuous second derivative, and $\int_{\bbx\in\ccalX} \phi''(\bbx) d\bbx < \infty$. \label{as:testfunction2}
\end{enumerate}
\end{assumption}
%
%%%%%%%%%%%%%%%%%%%%%%%%%%%%%%%%%%%%%%%%%%%%%%%%%%%%%%%%%%%%%%%%%%%%%%%%%%%%%%%%%%%%%%%%%%%%%%%%%%%%%%%%%%%%%%%%%%%%%%%%%%%%%%%%%%%%%%%%%%%
\begin{assumption}\label{as:2}
The kernel function associated with RKHS is such that  $\int_{\bbx \in \ccalX} \kappa_{\bbx^{(n)}}(\bbx) = 1$,  $\int_{\bbx \in \ccalX} \bbx \kappa_{\bbx^{(n)}} (\bbx) = 0$, and $\sigma_{\kappa}^2 = \int_{\bbx \in \ccalX} \bbx^2 \kappa_{\bbx^{(n)}}(\bbx) >0$.
\end{assumption}

\begin{assumption}\label{as:approximation}
Let $I_{N}(\phi)$ and $\hat I_{N}(\phi)$ be the integral estimators for test function $\phi$ associated with the uncompressed and compressed posterior densities. We define the approximation error for $\phi\notin\mathcal{F}$ where  $\mathcal{F}:=\{f\in\mathcal{H}~|~ \|f\|_{\mathcal{H}}\leq 1\}$, as 
 	 \begin{align}\label{approximation}
\!	I(\phi,f\!)\!= \!\sup_{|\phi| \leq 1}\!\! \left(\mathbb{E}[\hat I_N(\phi) \!-\!  {I_N(\phi)}]\! \right)\!-\! \sup_{f \in \mathcal{F} }\!\!\left(\mathbb{E}[\hat I_N(f) \!-\! {I_N(f)}] \right),
\end{align}
We assume that $I(\phi,f)\leq G$, where $G$ is a finite constant. 
\end{assumption}
%\red{skeptical about the above assumption}
%%%%%%%%%%%%%%%%%%%%%%%%%%%%%%%%%%%%%%%%%%%%%%%%%%%%%%%%%%%%%%%%%%%%%%%%%%%%%%%%%%%%%%%%%%%%%%%%%%%%%%%%%%%%%%%%%%%%%%%%%%%%%%%%%%%%%%%%%%%
%\begin{assumption}\label{as:last}
%\red{Denote as $\hat{m}$ and $\tilde{m}$ the mean embeddings (which belongs to $\ccalH$) of distributions $\hat{\rho},\tilde{\rho}$ \cite{sriperumbudur2010hilbert}} given by
%%
%\red{$$ \hat{m}=\mathbb{E}_{\hat{\rho}}[ \kappa_{\bbx}(\cdot) ] \; , \qquad \tilde{m} = \mathbb{E}_{\tilde{\rho}}[ \kappa_{\bby}(\cdot) ].$$}
%%
%%
%\red{The {'dist'} distance between the full distributions lower-bounds the distance between their mean embeddings: $dist(\hat{\rho},\tilde{\rho}) \leq  \| \hat{m} - \tilde m \|_{\ccalH}$, which are related by a multiplicative factor $dist(\hat{\rho},\tilde{\rho})  = K \| \hat{m} - \tilde m \|_{\ccalH}$.}
%%
%\end{assumption}
%
%
%\blue{need some comments about assumptions}
Assumption \ref{as:first}\ref{as:testfunction1} is a textbook condition in the analysis of Monte Carlo methods, and appears in \cite{mcbook}. Assumptions \ref{as:first}\ref{as:testfunction2} and \ref{as:2} are required conditions for establishing consistency of kernel density estimates and are standard -- see \cite{wasserman2006all}[Theorem 6.28]. 
%Assumption \ref{as:last} is {reasonable} considering the fact that a distribution is completely characterized by its moments \cite{durrett2019probability}. Hence if the full distribution is close according to some metric, then the moments are close -- see \cite{scholkopf2015computing}[eqn. (14)]. This assumption formalizes this statement for  means and metrics in RKHSs \cite{sriperumbudur2010hilbert}.
%
 We begin by noting that under Assumption \ref{as:first}, we have classical statistical consistency of importance sampling as the number of particles becomes large as stated in Lemma \ref{first_theorem} in Appendix \ref{first_proof}.
%
%The proof is in Appendix \ref{first_proof}. 
%
This result enables characterizing the bias of Algorithm \ref{alg:compressed_online_kernel_importance_sampling}, given next in Lemma \ref{intermidiate_lemma}. Assumption \ref{as:approximation} is non-standard and we use it to bound the error due to continuity conditions imposed by operating in the RKHS. We note that if $\phi\in\mathcal{F}$ which is the case for most practical applications, then $G=0$. 
%\blue{we need to state some assumptions here}
%
%

\begin{lemma}\label{intermidiate_lemma}
Define the second moment of the true unnormalized density $\rho$ as in Lemma \ref{first_theorem}. Then, under Assumptions \ref{as:first}-\ref{as:approximation}, the estimator of Alg. \ref{alg:compressed_online_kernel_importance_sampling} satisfies
\begin{align}\label{lemma}
	\big|\sup_{|\phi| \leq 1}  \big(\mathbb{E}[\hat I_N(\phi)- I(\phi)] \big) \big| \leq  \sum_{n=1}^{N}\epsilon_n+\frac{24}{N} \rho +G.
\end{align}
where $\epsilon_n$ is the compression budget for each $n$.
\end{lemma}

\begin{myproof}
	Inspired by \cite{agapiou2017importance}, begin by denoting $\hat I_N(\phi)$ as the integral estimate given by Algorithm \ref{alg:compressed_online_kernel_importance_sampling}. Consider the bias of the integral estimate ${\hat I_N(\phi) - I(\phi)}$, and add and subtract $I_N(\phi)$, the uncompressed normalized importance estimator that is the result of Algorithm \ref{alg:importance_sampling}, to obtain
	\begin{align}\label{new}
		\hat I_N(\phi) - I(\phi)=& \hat I_N(\phi) -I_N(\phi)+I_N(\phi) - I(\phi).
	\end{align}
	%
	%where $\tilde I_N(\phi)$ is the uncompressed integral estimation at $N$ (\red{We might need to define it properly because without intorducing the kernalized IS, we have directly introduced the compressed kernelized density in Algorithm \ref{alg:compressed_online_kernel_importance_sampling}.}) $I_N(\phi)$ is the integral estimation given by Algorithm \ref{alg:importance_sampling}.  
	%
	Take the expectation on both sides with respect to the population posterior \eqref{eq:posterior} to obtain
	\begin{align}\label{new2}
			\mathbb{E}[\hat I_N(\phi) - I(\phi)]=& \mathbb{E}[\hat I_N(\phi) -  I_N(\phi)]+\mathbb{E}[I_N(\phi) - I(\phi)].
	\end{align}
	%
	%where we add and subtract $\mathbb{E}[I(\phi)]$ inside the second-to-last term and group like terms.
	%
	Let's compute the $\sup$ of both sides of \eqref{new2} over range $|\phi | \leq 1$ and use the fact that a sup of a sum is upper-bounded by the sum of individual terms:
	%
	%Using the triangle inequality, we can write 
	%
	\begin{align}\label{new3}
		\sup_{|\phi | \leq 1} \left(\mathbb{E}[\hat I_N(\phi) - I(\phi)] \right)&\leq  \sup_{|\phi | \leq 1} \left(\mathbb{E}[\hat I_N(\phi) -  I_N(\phi)] \right)\nonumber
		\\
		%& + \sup_{|\phi | \leq 1} \left(\mathbb{E}[\tilde I_N(\phi)-I(\phi)]\right) \nonumber \\
		&\quad+\! \sup_{|\phi| \leq 1}\!\!\left([\mathbb{E}[I_N(\phi) \!-\! I(\phi\!)] \right)\!.
	\end{align}
Now add and subtract the supremum over the space $\mathcal{F}:=\{f\in\mathcal{H}~|~ \|f\|_{\mathcal{H}}\leq 1\}$ 
		to the first term on the right hand side of \eqref{new3} to write
		\begin{align}\label{IPM0_0}
			&\sup_{|\phi| \leq 1} \Big(\mathbb{E}[\hat I_N(\phi) \!-\! {I(\phi)}] \Big) \\
		&	\leq \sup_{f \in \mathcal{F} }\left(\mathbb{E}[\hat I_N(f) \!-\! {I_N(f)}] \right)+ \sup_{|\phi| \leq 1}\!\!\left([\mathbb{E}[I_N(\phi) - I(\phi)] \right)\nonumber \\
		&+\underbrace{\sup_{|\phi| \leq 1} \left(\mathbb{E}[\hat I_N(\phi) \!-\!  {I_N(\phi)}]\! \right)- \sup_{f \in \mathcal{F} }\left(\mathbb{E}[\hat I_N(f) \!-\! {I_N(f)}] \right)}_{I(\phi,f)\leq G} \nonumber 
		\end{align}
Observe that the last line on the right-hand side of the preceding expression defines the integral function approximation error $I(\phi,f)$ defined in \eqref{approximation}, which is upper-bounded by constant $G$ (Assumption \ref{as:approximation}).
Now, compute the absolute value of both sides of \eqref{IPM0_0}, and to the second term on the first line, pull the absolute value inside the supremum.
	%
	%$$|\sup \left(\mathbb{E}[\tilde I_N(\phi)-I_N(\phi)]\right)  | \leq  \sup \left|\mathbb{E}[\tilde I_N(\phi)-I_N(\phi)]\right|.$$
	Doing so allows us to apply \eqref{eq:posterior_contraction} (Lemma \ref{first_theorem}) to the second term on the first line, the result of which is:
	\begin{align}\label{new4}
		\big|\sup_{|\phi| \leq 1}  \big(\mathbb{E}[\hat I_N(\phi)- I(\phi)] \big) \big| \leq & \sup_{f \in \mathcal{F} }\left(\mathbb{E}[\hat I_N(f) \!-\! {I_N(f)}] \right)\nonumber
		\\
		&+\frac{24}{N} \rho+G,
	\end{align}
	where $G$ is defined in Assumption \ref{as:approximation} and note that $G=0$ if $\phi \in \mathcal{H}$. 
		It remains to address the first term on the right hand side of \eqref{new4}, which noticeably defines an instance of an integral probability metric (IPM) [Muller 1997], i.e., 
		\begin{align}\label{IPM1}
			\sup_{f \in \mathcal{F} } &\left(\mathbb{E}[\hat I_N(f) \!-\! {I_N(f)}] \right)\nonumber
			\\
			&=\sup_{f \in \mathcal{F} } \left(\int f(\bbx)d\hat \mu_N -\int f(\bbx)d\tilde\mu_N\right),
	\end{align} 
	where $\hat\mu_N$  is the pre-image unnormalized density estimate obtained by solving \eqref{eq:proximal_hilbert_representer}. The IPM in \eqref{IPM1} is exactly equal to Maximum Mean Discrepancy (MMD)\cite{fortet1953convergence} for test functions  in the RKHS $f\in\mathcal{H}$. This observation allows us to write 
	\begin{align}\label{eq:mmd}
			\sup_{f \in \mathcal{F} } \left(\mathbb{E}[\hat I_N(f) \!-\! {I_N(f)}] \right) =\|\beta_{ N}-\gamma_{N}\|_{\mathcal{H}}.
	\end{align}
	{where $\gamma_{N}$ is the kernel mean embedding corresponding to the uncompressed measure estimate $ \tilde\mu_N$. Next, consider the term $\|\beta_{ N}-\gamma_{N}\|_{\mathcal{H}}$ and add and subtract the per step uncompressed estimate $\tilde\beta_{N}$ [cf. \eqref{eq:importance_sampling_stack2}] to obtain
	}
	\begin{align}\label{triangle}
			\|\beta_{N}-\gamma_{N}\|_{\mathcal{H}} =& \|(\beta_{ N}-\tilde\beta_{ N})+(\tilde\beta_{ N}-\gamma_{N})\|_{\mathcal{H}}\nonumber
			\\
			\leq & \|\beta_{ N}-\tilde\beta_{N}\|_{\mathcal{H}}+\|\tilde\beta_{ N}-\gamma_{\nu_N}\|_{\mathcal{H}}\nonumber
			\\
			\leq &  \epsilon_N+\|\tilde\beta_{N}-\gamma_{N}\|_{\mathcal{H}},
	\end{align}
	where the last inequality holds from the fact that $\|\beta_{N}-\tilde\beta_{ N}\|_{\mathcal{H}}\leq \epsilon_N$ (cf. Algorithm \ref{alg:komp}). Now we substitute the values of $\tilde\beta_{N}$ and $\gamma_{N}$ using the update in \eqref{eq:kernel_mean_emb2}, and we get 
	\begin{align}\label{triangle2}
			\|\beta_{ N}-\gamma_{N}\|_{\mathcal{H}}
			\leq &  \epsilon_N+\|\beta_{ {N-1}}-\gamma_{{N-1}}\|_{\mathcal{H}}.
	\end{align}
	Using the above recursion, we easily obtain
	{\begin{align}\label{triangle3}
			\|\beta_{ N}-\gamma_{N}\|_{\mathcal{H}}
			\leq &  \sum_{n=1}^{N}\epsilon_n.
		\end{align}
Hence, using the upper bound of \eqref{triangle3} in \eqref{eq:mmd}, and then substituting the result into \eqref{new4} yields 
		\begin{align}\label{new5}
			\big|\sup_{|\phi| \leq 1}  \big(\mathbb{E}[\hat I_N(\phi)- I(\phi)] \big) \big| \leq  \sum_{n=1}^{N}\epsilon_n+\frac{24}{N} \rho +G.
		\end{align}}
		as stated in Lemma \ref{intermidiate_lemma}.
\end{myproof}

{With this technical lemma in place, we are ready to state the main result of this paper. }
    
%
%Observe, however, that for consistency \eqref{eq:consistent}, we require that the number of generated particles from the importance distribution and the parameterization of the importance distribution to be equal $\{\bby_u\}$, i.e. $M_n=n$. Therefore, for streaming settings where the number of observations is to be infinite, i.e., $t\rightarrow \infty$, we require infinitely many samples $\{\bbx^{n} \}$ from the importance distribution, whose parameterization in terms of samples $\bbx^{n}$ becomes infinite. The question, then, is how do we obtain estimators for \eqref{eq:template_problem} which are \emph{nearly} consistent but of finite memory.
%%%%%%%%%%%%%%%%%%%%%%%%%%%%%%%%%%%%%%%%%%%%%%%%%%%%%%%%%%%%%%%%%%%%%%%%%%%%%%%%%%%%%%%%%%%%%%%%%%%%%%%%%%%%%%%%%%%%%%%%%%%%%%%%%%%%%%%%%%%%%%%%%%%%%%% %%%%%%%%%%%%%%%%%%T  H  E  O  R  E  M%%%%%%%%%%%%%%%%%%%%%%%%%%%%%%%%%%%%%%%%%%%%%%%%%%%%%%%%%%%%%%%%%%%%%%%%%%%%%%%%%%%%%%
%%%%%%%%%%%%%%%%%%%%%%%%%%%%%%%%%%%%%%%%%%%%%%%%%%%%%%%%%%%%%%%%%%
%\subsection{Theoretical guarantees for Compressed Internalized IS (CKIS) }
\begin{theorem}\label{thm:main_result}
%
%  Suppose $\pi$, the proposal distribution is absolutely continuous with respect to $q$, the population posterior, and both are probability measures defined over $\ccalX$. Then they are related by their Radon-Nikodyn derivative:
%$$
%\frac{d q}{d \pi}(x) : = \frac{g(x ) }{\int_{\ccalX} g(x) \pi(dx)}
%$$
%where $g$ is the unnormalized density (Radon-Nikodyn derivative) of $q$ with respect to $\pi$. Define the second moment of the Radon-Nikodyn derivative as 
%%
%$$\rho:= \frac{\pi(g^2)}{\mu(g^2)}$$
%%
Define $\rho = \frac{\pi(g^2)}{q(g^2)}$ as the variance of the unnormalized importance density with respect to importance weights $g$ as  in Lemma \ref{first_theorem}. Then under Assumptions \ref{as:first}-\ref{as:approximation}, we have the following approximate consistency results: 
\begin{enumerate}[label=(\roman*)] 
\item \label{statement:diminish} for diminishing compression budget $\epsilon_n=\alpha^n$ with $\alpha\in(0,1)$, the estimator of Alg. \ref{alg:compressed_online_kernel_importance_sampling} satisfies
\begin{align}\label{thm:diminish}
	\big|\sup_{|\phi| \leq 1}  \big(\mathbb{E}[\hat I_N(\phi)- I(\phi)] \big) \big| \leq  \frac{\alpha}{1-\alpha}+\mathcal{O}\left(\frac{1}{N}\right) +G.
\end{align}
To obtain a $\delta$ accurate integral estimate, we need at least $N\geq\mathcal{O}\left(\frac{1}{\delta}\right)$ particles and compression attenuation rate sufficiently large such that $0<\alpha\leq 1/(1+(2/\delta))$.
\item  \label{statement:constant} for constant compression budget $\epsilon_n=\epsilon>0$ and memory $\mathcal{M}:=\mathcal{O}\left(\frac{1}{\epsilon^{1/(2p)}}\right)$
\begin{align}\label{final0}
	\big|\sup_{|\phi| \leq 1}  \big(\mathbb{E}[\hat I_N(\phi)- I(\phi)] \big) \big| \leq \mathcal{O}\left(\frac{N}{\mathcal{M}^{1/2p}}+\frac{1}{N}\right)
\end{align}
and
\begin{align}
	\mathcal{O}\left(\frac{1}{\delta}\right)\leq N\leq \mathcal{O}\left(\delta {\mathcal{M}^{1/(2p)}}\right),
\end{align}
which implies that $\mathcal{M}\geq\mathcal{O}\left(\frac{1}{\delta^{4p}}\right)$. 
\end{enumerate}
\end{theorem}
\begin{myproof}
Consider the statement of Lemma \ref{lemma} and to proceed next, we characterize the behavior of the term $\sum_{n=1}^{N}\epsilon_n$ since it eventually determines the final bias in the integral estimation. 

{\bf \noindent Theorem \ref{thm:main_result}\ref{statement:diminish}}: Diminishing compression budget: Let us consider $\epsilon_n=(\alpha)^n$ with $\alpha=(0,1)$, which implies that
\begin{align}\label{diminishing}
\sum_{n=1}^{N}\epsilon_n= \frac{\alpha(1-\alpha^N)}{1-\alpha}\leq \frac{\alpha}{1-\alpha}.
\end{align}
Substituting \eqref{diminishing} into the right hand side of \eqref{new5}, we get
\begin{align}\label{new6}
	\big|\sup_{|\phi| \leq 1}  \big(\mathbb{E}[\hat I_N(\phi)- I(\phi)] \big) \big| \leq  \frac{\alpha}{1-\alpha}+\frac{24}{N} \rho.
\end{align}
To obtain a $\delta$ accurate integral estimate, we need $N\geq\frac{48\rho}{\delta}$ and $0<\alpha\leq 1/(1+(2/\delta))$.

 \noindent{ \bf{Theorem \ref{thm:main_result}\ref{statement:constant}}}: Constant compression budget: Let us consider $\epsilon_n=\epsilon$, which implies that
	\begin{align}\label{constant_epsilon}
		\sum_{n=1}^{N}\epsilon_n= N\epsilon.
	\end{align}
Using \eqref{constant_epsilon} in \eqref{new5}, we can write
\begin{align}\label{final}
	\big|\sup_{|\phi| \leq 1}  \big(\mathbb{E}[\hat I_N(\phi)- I(\phi)] \big) \big| \leq  N\epsilon+\frac{24}{N} \rho.
\end{align}
For a constant compression budget $\epsilon$, from Theorem \ref{thm:model_order}, we have 
\begin{align}
		M_{\infty}\leq \mathcal{O}\left(\frac{1}{\eps^{2p}}\right):=\mathcal{M}.
\end{align}
If we are given a maximum memory requirement $\mathcal{M}$, then we can choose $\epsilon$ as  
\begin{align}
	\epsilon=\frac{G}{\mathcal{M}^{1/(2p)}},
\end{align}
where $G$ is a bound on the unnomalized weight $g(x(u))\leq G$ for all $u$. Using this lower bound value of $\epsilon$ in \eqref{final}, we get
\begin{align}\label{final1}
	\big|\sup_{|\phi| \leq 1}  \big(\mathbb{E}[\hat I_N(\phi)- I(\phi)] \big) \big| \leq \frac{NG}{\mathcal{M}^{1/2p}}+\frac{24}{N} \rho.
\end{align}
%%%
\begin{figure}[h] 
	\centering 
	\includegraphics[scale=0.15]{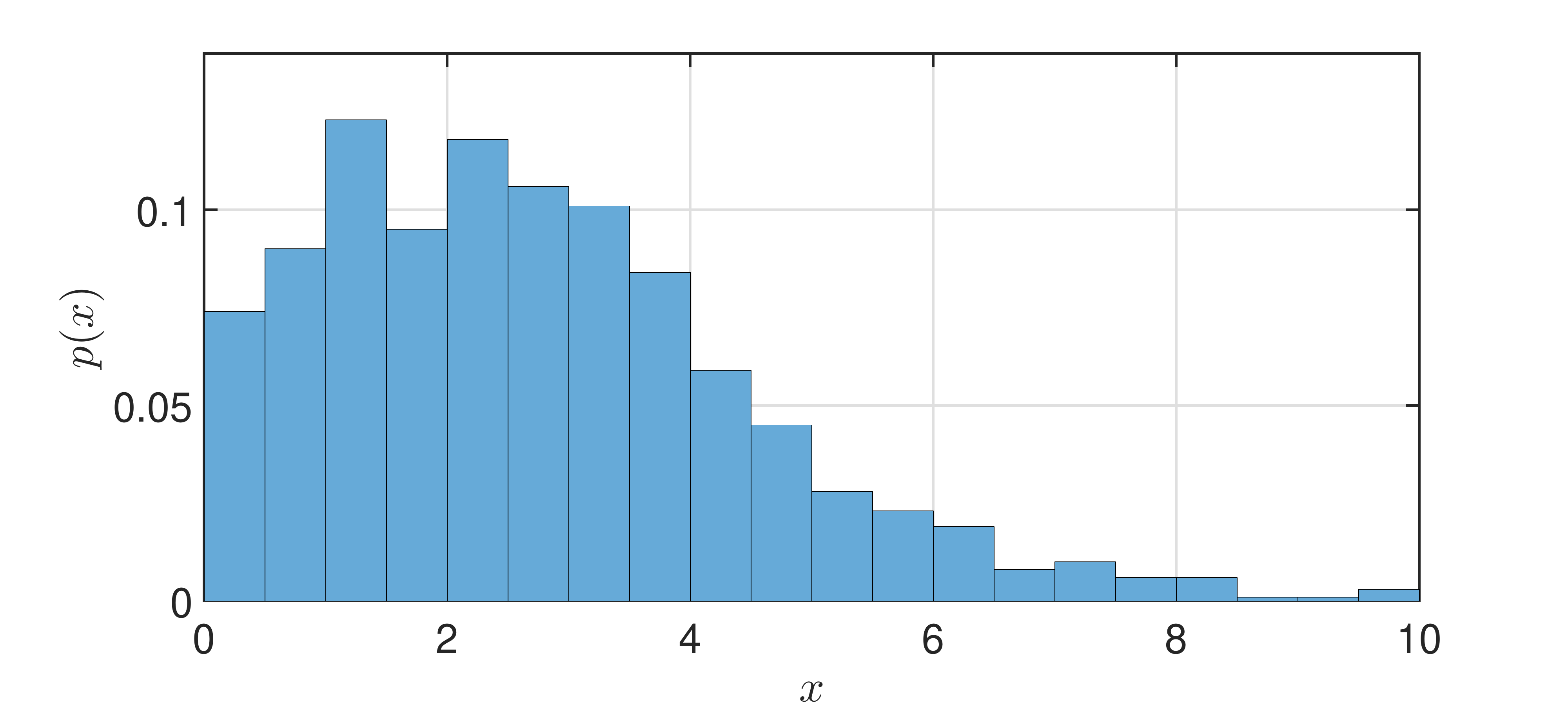}
	\caption{Particle histogram for direct sampling experiment. } 
	\label{hist}
	\vspace{-8mm}
\end{figure}
%  

%
%%%%%%%%%%%%%%%%%%%%%%%%%%%%%%%%%%%%%%%%%%%%%%%%%%%%%%%%%%%%%%%%%%%%%%%%%%%%%%%%%%%%%%%%%%%%%%%%%%%%%%%%%%%%%%%%%%%%%%%%%%%%%%%%%%%%%%%%%%%%%%%%%%%%%%%%%%%%%%%%%%%%%%%%%%%%%%%%%%%%%%%%%%%%%%%%%%%%%%%%%%%%%%%%%%
  \begin{figure*}[t]
	\centering
%	\begin{subfigure}[b]{0.475\textwidth}
%		\centering
%		\includegraphics[width=\textwidth]{figs/first.eps}
%		\caption{{\small Unnormalized intergral estimate vs. particle $n$\vspace{-2mm}}}    
%		\label{fig:unnormalized}
%	\end{subfigure}
%	\hfill
	\begin{subfigure}[b]{0.32\textwidth}  
		\centering 
		\includegraphics[width=\textwidth]{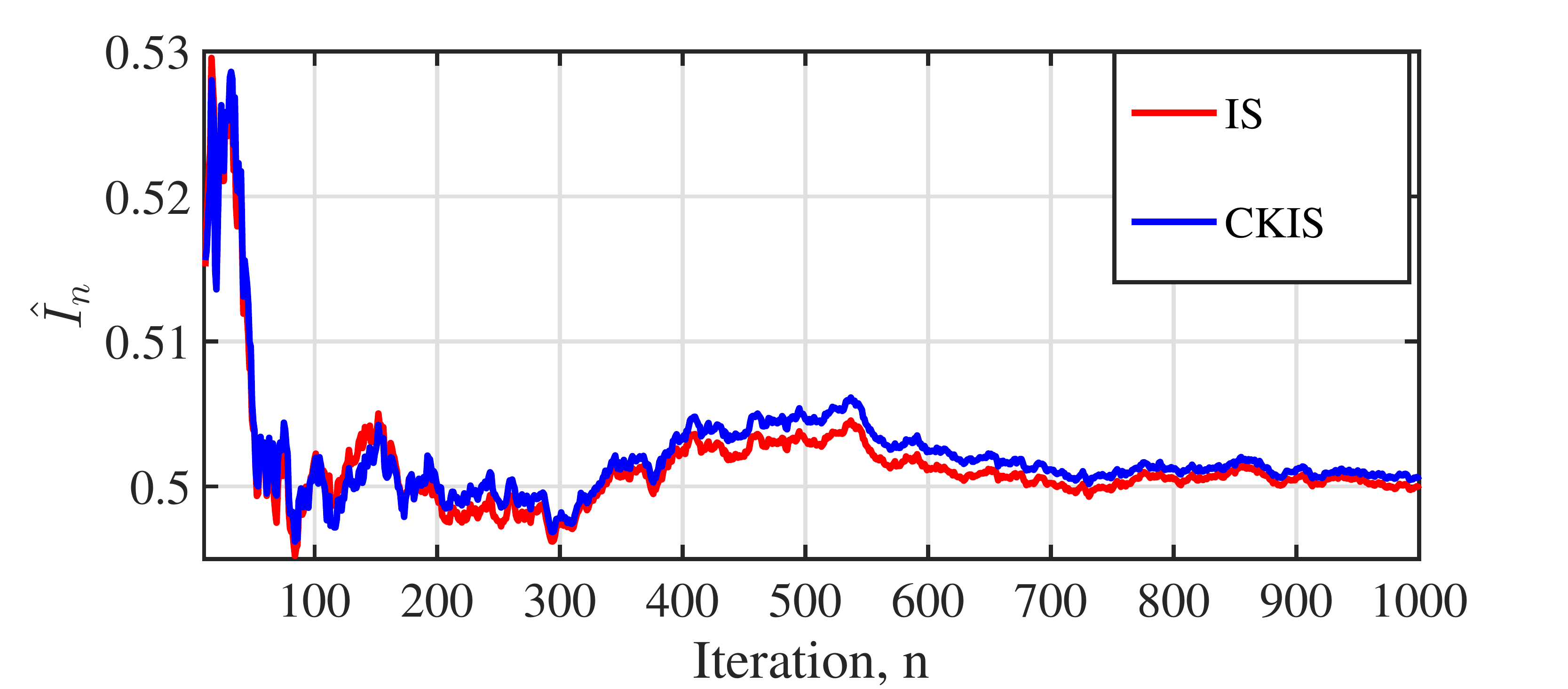}
		\caption[\red{y label}]%
		{{\small Integral estimate. \vspace{-0mm}}}    
		\label{fig:normalized}
	\end{subfigure}
	\begin{subfigure}[b]{0.32\textwidth}   
		\centering 
		\includegraphics[width=\textwidth]{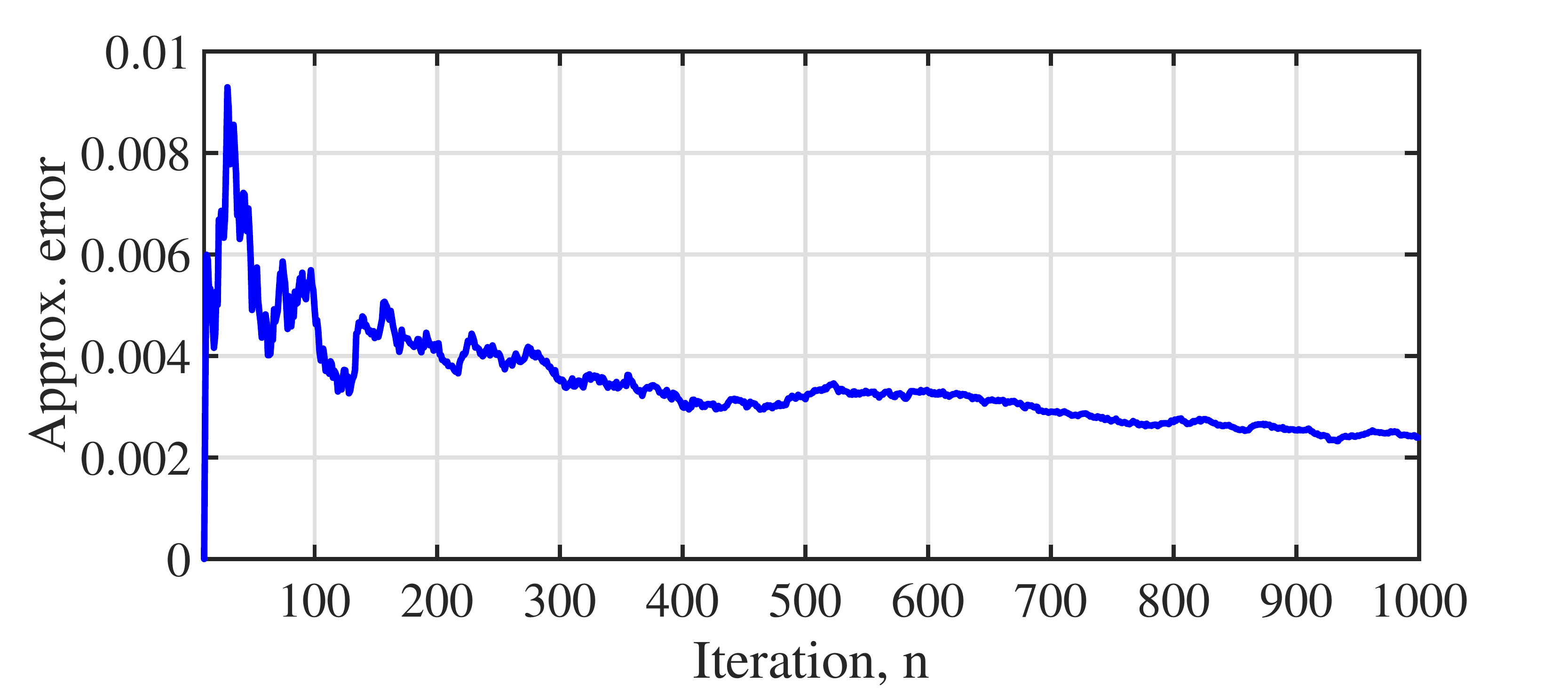}
		\caption[]%
		{{\small Integral Error between Algs. \ref{alg:importance_sampling} - \ref{alg:compressed_online_kernel_importance_sampling}. }}    
		\label{fig:error}
	\end{subfigure}
	\begin{subfigure}[b]{0.32\textwidth}   
		\centering 
		\includegraphics[width=\textwidth]{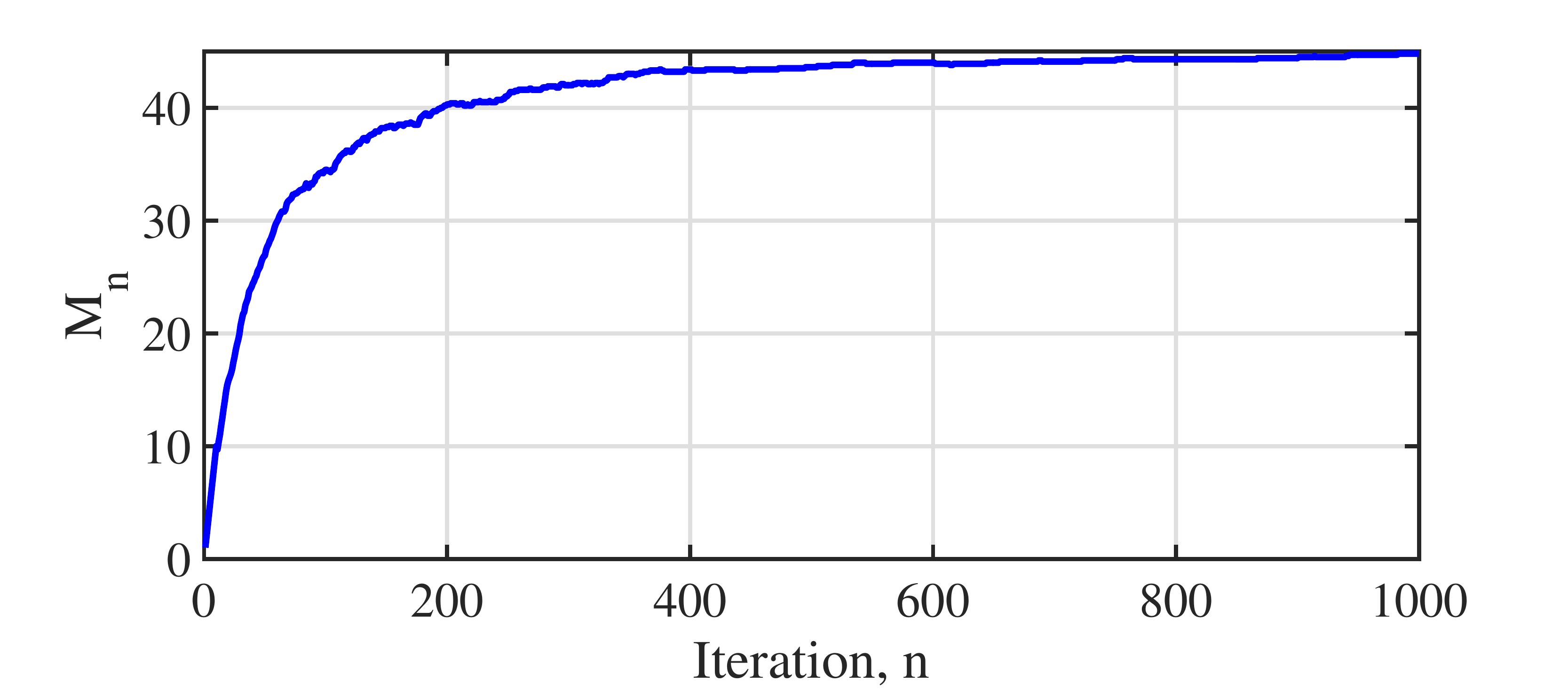}
		\caption[]%
		{{\small Model Order. }}    
		\label{fig:model_order}
	\end{subfigure}
	\caption[ The average and standard deviation of critical parameters ]
	{\small Simulation results for Alg. \ref{alg:compressed_online_kernel_importance_sampling} run with Gaussian kernel ($h=0.01$) and compression budget $\epsilon=3.5$ for the problem \eqref{eq:sample_problem2}. The memory-reduction scheme nearly preserves statistical consistency, while yielding reasonable complexity, whereas Alg. \ref{alg:importance_sampling} attains exact consistency as its memory grows unbounded with index $n$. {All the plots are averaged over $10$ iterations.}} 
	\label{fig1}
	\vspace{-0mm}
\end{figure*}

Note that the first term in the above expression increases with $N$ and the second term decreases with $N$, hence we obtain a tradeoff between memory and accuracy for the importance sampling based estimator. For a given memory $\mathcal{M}$ (number of elements we could store in the dictionary) and the required accuracy $\delta$, we obtain the following bound on the number of iterations $N$ 
\begin{align}
	\frac{48\rho}{\delta}\leq N\leq \frac{\delta {M^{1/(2p)}}}{2G},
\end{align}
which implies that for increased accuracy we need to run for more iterations but need more memory,  and vice versa.   
\end{myproof}
%
%$\hfill \blacksquare$

Theorem \ref{thm:main_result}  establishes that the compressed kernelized importance sampling scheme proposed in Section \ref{sec:algorithm} is \emph{nearly} %\green{This is the first time the word "nearly" appear since the title and abstract (besides the caption above!). I think it would deserve a clear explanation in the intro and also earlier in this section.}
asymptotically consistent. Note that the right hand side in \eqref{thm:diminish} consists of three terms depending upon $\alpha$, $\frac{1}{N}$, and $G$. If we ignore $G$, which actually depends upon the approximation associated with function $\phi(\cdot)$, the other two terms can be made arbitrarily small by making $\alpha$ close to zero and a very high $N$. Hence, the integral estimation can be made arbitrarily close to exact integral. However, when these parameters are fixed positive constants, they provide a tunable tradeoff between bias and memory. That is, when the compression budget is a positive constant, then the memory of the posterior distribution representation is finite, as we formalize next. 
%%%%%%%%%%%%%%%%%%%%%%%%%%%%%%%%%%%%%%%%%%%%%%%%%%%%%%%%%%%%%%%%%%%%%%%%%%%%%%%%%%%%%%%%%%%%%%%%%%%%%%%%%%%%%%%%%%%%%%%%%%%%%%%%%%%%%%%%%%%%%%%%%%%%%%% %%%%%%%%%%%%%%%%%%T  H  E  O  R  E  M%%%%%%%%%%%%%%%%%%%%%%%%%%%%%%%%%%%%%%%%%%%%%%%%%%%%%%%%%%%%%%%%%%%%%%%%%%%%%%%%%%%%%%
%%%%%%%%%%%%%%%%%%%%%%%%%%%%%%%%%%%%%%%%%%%%%%%%%%%%%%%%%%%%%%%%%%
%  
\begin{theorem}\label{thm:model_order}
 Under Assumptions \ref{as:first}-\ref{as:2} (in Section \ref{apx_technicalities}), for compact feature space $\ccalX$ and bounded importance weights $g(\bbx^{(n)})$, the model order $M_n$ for Algorithm \ref{alg:compressed_online_kernel_importance_sampling}, for all $n$ is bounded by 
\begin{align}
1\leq 	M_{n}\leq \mathcal{O}\left(\frac{1}{\eps^{2p}}\right).\end{align}
\end{theorem}
Theorem \ref{thm:model_order} (proof in Appendix \ref{third_proof}) contrasts with the classical bottleneck in the number of particles required to represent an arbitrary posterior, which grows unbounded \cite{agapiou2017importance}. While this is a pessimistic estimate, experimentally we observe orders of magnitude reduction in complexity compared to exact importance sampling, which is the focus of the subsequent section. Observe, however, that this memory-reduction does not come for free, as once the compression budget is fixed, the memory is fixed by the ratio $\frac{1}{\epsilon^{2p}}$ that eventually results in a lower bound on the accuracy of the integral estimate.
%%%%%%%%%%%%%%%%%%%%%%%%%%%%%%%%%%%%%%%%%%%%%%%%%%%%%%%%%%%%%%%%%%%%%%%%%%%%%%%%%%%%%%%%%%%%%%%%%%%%%%%%%%%%%%%%%%%%%%%%%%%%%%%%%%%%%%%%%%%%%%%%%%%%%%%%%%%%%%%%%%%%%%%%%%%%%%%%%%%%%%%%%%%%%%%%%%%%%%%%%%%%%%%%%%%%%%%%%%%%%%%%%%%%%%%%%%%%%%%%%%%%
%%%%%%%%%%%%%%%%%%%%%%%%%%%%%%%%%%%%%%%%%%%%%%%%%%%%%%%%%%%%%%%%%%%%%%%%%%%%%%%%%%%%%%%%%%%%%%%%%%%%%%%%%%%%%%%%%%%%%%%%%%
 \begin{figure*}[t]
        \centering
        \begin{subfigure}[b]{0.32\textwidth}  
            \centering 
            \includegraphics[width=\textwidth]{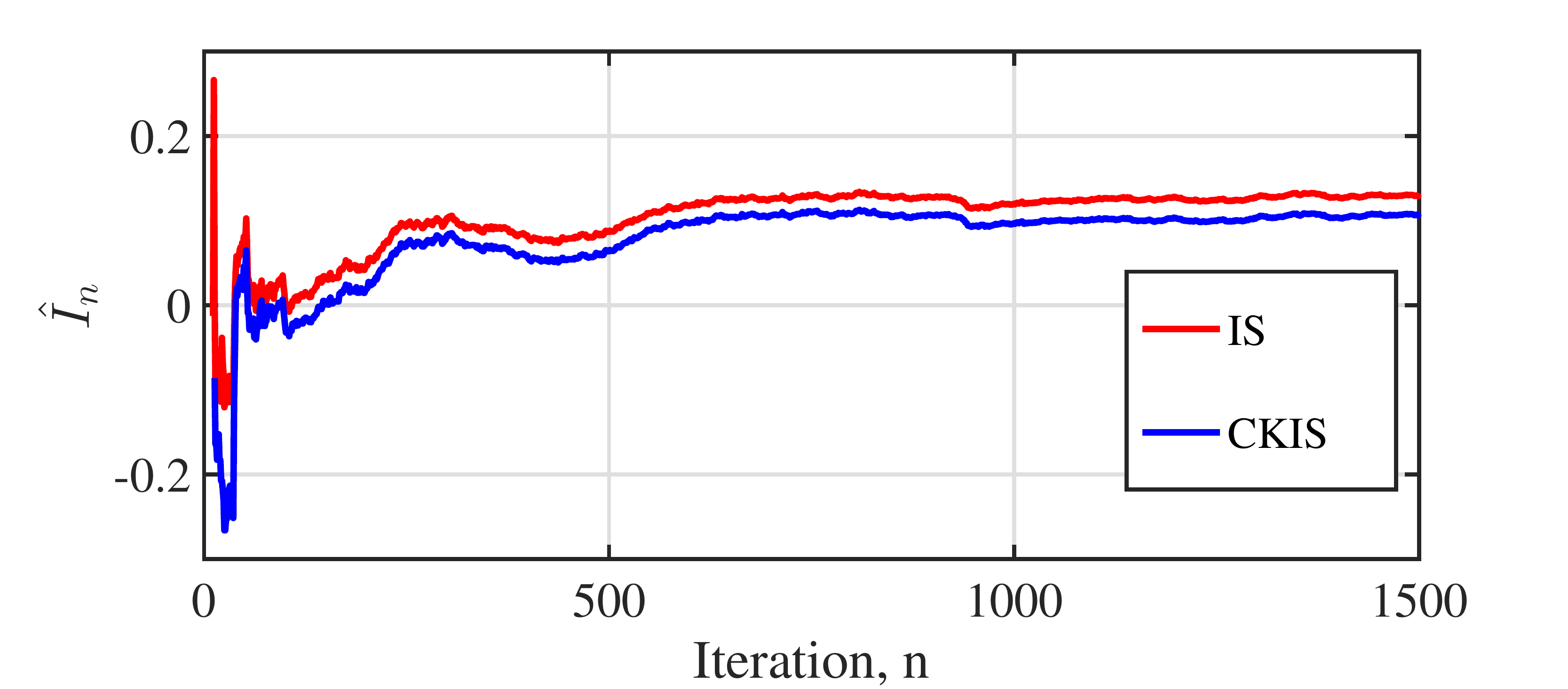}
            \caption[]%
            {{\small Integral estimate. \vspace{-0mm}}}    
            \label{fig:normalized1}
        \end{subfigure}
        \begin{subfigure}[b]{0.32\textwidth}   
            \centering 
            \includegraphics[width=\textwidth]{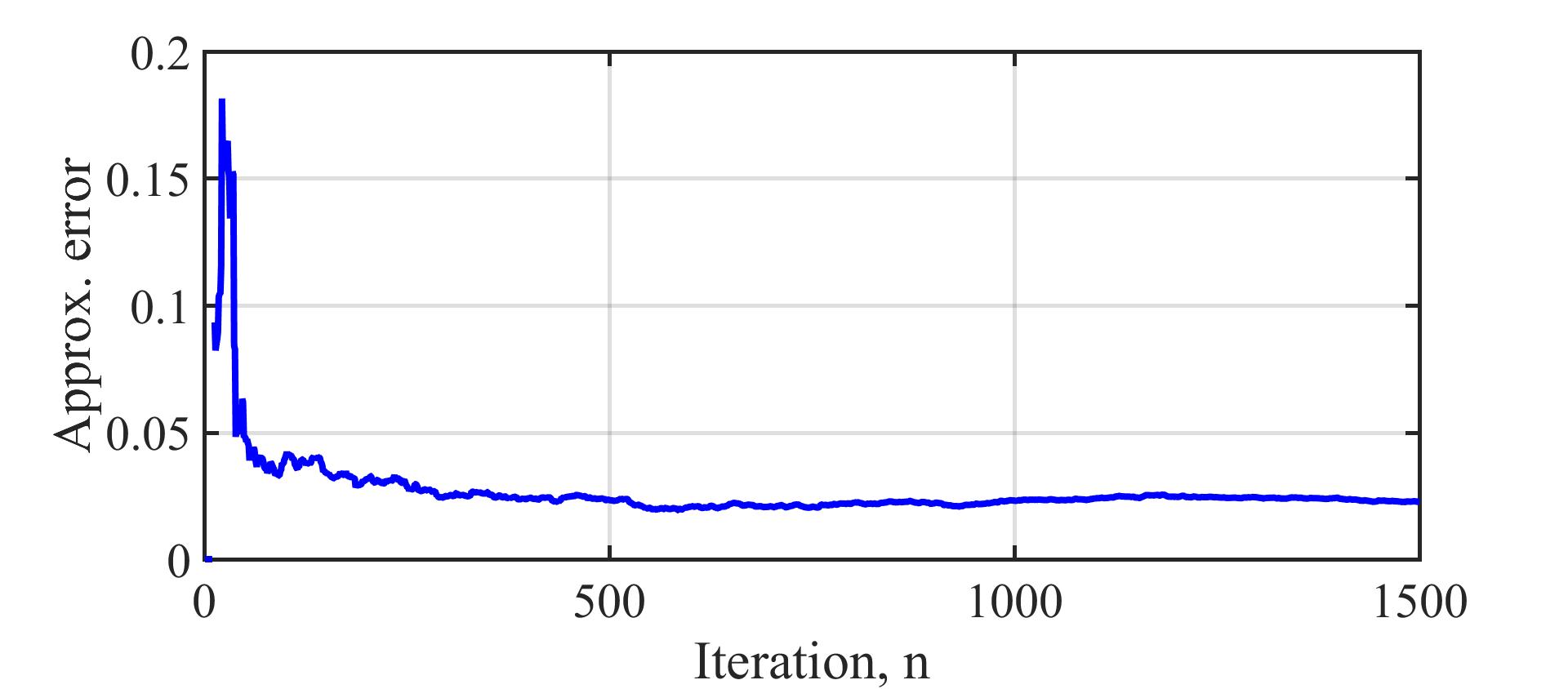}
            \caption[]%
            {{\small Integral Error between Algs. \ref{alg:importance_sampling} - \ref{alg:compressed_online_kernel_importance_sampling}. }}    
            \label{fig:error1}
        \end{subfigure}
        \begin{subfigure}[b]{0.32\textwidth}   
            \centering 
            \includegraphics[width=\textwidth]{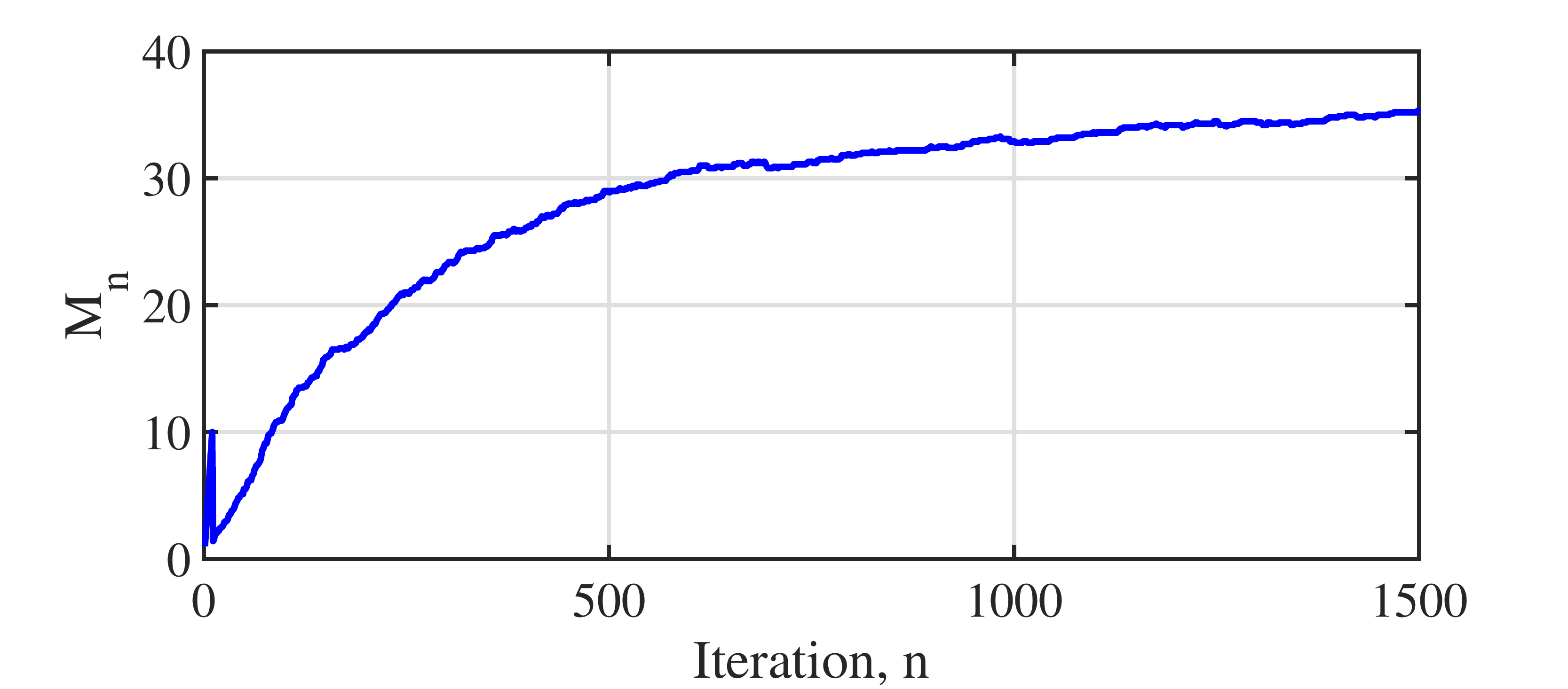}
            \caption[]%
            {{\small Model Order. }}    
            \label{fig:model_order1}
        \end{subfigure}
        \caption[ The average and standard deviation of critical parameters ]
        {\small Simulation results for Alg. \ref{alg:compressed_online_kernel_importance_sampling} with indirect IS,  run with Gaussian kernel ($h=0.012$) and compression budget $\epsilon=10^{-3}$ for the problem \eqref{eq:sample_problem2}. The memory-reduction scheme nearly preserves statistical consistency, while yielding reasonable complexity, whereas Alg. \ref{alg:importance_sampling} attains exact consistency as its memory grows unbounded with index $n$. {All the plots are averaged over $10$ iterations.}\vspace{-4mm}} 
        \label{fig2}
    \end{figure*}
    
\section{Experiments}\label{sec:experiments}
 
 \subsection{Direct Importance Sampling}\label{direct_IS}
In this section, we conduct a simple numerical experiment to demonstrate the efficacy of the proposed algorithm in terms of balancing model parsimony and statistical consistency. We consider the problem of estimating the expected value of function $\phi(\bbx)$
%\begin{align}\label{eq:sample_problem}
%\phi(x)=2\sin\left(\frac{\pi}{(1.5x)}\right)
%\end{align}
with the target $q(\x)$ and the proposal $\pi(\x)$ given by 
\begin{align}\label{eq:sample_problem2}
\!\!\phi(x)=2\sin\!\left(\frac{\pi}{(1.5x)}\!\right)\; , \ q(\x)\!=\!\frac{1}{\sqrt{2\pi}}\exp\left(\!\!-\frac{(x-1)^2}{2}\right)\; , \ \nonumber\\
\pi(\x)\!=\!\frac{1}{\sqrt{4\pi}}\exp\left(\!\!-\frac{(x-1)^2}{4}\right),
\end{align}
to demonstrate that generic Monte Carlo integration allows one to track generic quantities of random variables that are difficult to compute under more typical probabilistic hypotheses. For \eqref{eq:sample_problem2}, since $q(\x)$ is known, this is referred to as ``direct importance sampling". We run Algorithm \ref{alg:importance_sampling}, i.e., classic importance sampling, and Algorithm \ref{alg:compressed_online_kernel_importance_sampling} for the aforementioned problem. For Algorithm \ref{alg:compressed_online_kernel_importance_sampling}, we select compression budget $\epsilon=3$, and used a Gaussian kernel with bandwidth $h=0.01$.  We track the normalized integral estimate  \eqref{eq:SNIS_estimator}, absolute integral approximation error, and the number of particles that parameterize the empirical measure (model order). 

We first represent the histogram of the particles generated in Fig.~\ref{hist}.  In Fig. \ref{fig:normalized}, we plot the un-normalized integral approximation error for Algorithms \ref{alg:importance_sampling} - \ref{alg:compressed_online_kernel_importance_sampling}, which are close, and the magnitude of the difference depends on the choice of compression budget. Very little error is incurred by kernel mean embedding and memory-reduction. The magnitude of the error relative to the number of particles generated is displayed in Fig. \ref{fig:error}: observe that the error settles on the order of $10^{-3}$. In Fig. \ref{fig:model_order}, we display the number of particles retained by Algorithm \ref{alg:compressed_online_kernel_importance_sampling}, which stabilizes to around $56$, whereas the complexity of the empirical measure given by Algorithm \ref{alg:importance_sampling} grows linearly with sample index $n$, which noticeably grows \emph{unbounded}.
% 

%%%%%%%%%%%%%%%%%%%%%%%%%%%%%%%%%%%%%%%%%%
 \subsection{Indirect Importance Sampling}\label{indirect_IS}
 As discussed in Sec. \ref{sec:prob}, in practice we do not know the target distribution $q(\x)$ and hence we use Bayes rule as described in \eqref{eq:posterior1} to calculate $q(\x^{(n)})$ at each instant $t$. We consider the observation model $y_t=b+\text{sin}(2\pi x)+\eta_t$ where $\eta_t\sim\mathcal{N}(0,\sigma^2)$. From the equality in \eqref{eq:posterior1}, we need the likelihood and a prior distribution to calculate $q(x^{(n)})$ using Bayes Rule [cf. \eqref{eq:posterior}]. Here we fix the likelihood (measurement model) and prior as
 %
 %\begin{align}
 $\mathbb{P}\left(\{y_k\}_{k\leq K} \given  x^{(n)}  \right)=\frac{1}{(2\pi\sigma_1^2)^{K/2}}\exp\left(-\frac{\|\y-x^{(n)}\|^2}{2\sigma_1^2}\right)$, $
  \mathbb{P}\left(x^{(n)}\right)=\frac{1}{(2\pi\sigma_2^2)}\exp\left(-\frac{(x^{(n)})^2}{2\sigma_2^2}\right). $
 %\end{align}
 We set $K=10$, $b=5$, $\sigma=0.1$, $\sigma_1=0.4$, $\sigma_2=1.6$, and compression budget $\epsilon=10^{-3}$. A uniform distribution $\mathcal{U}[3,7]$ is used as the importance distribution. The results are reported in Fig. \ref{fig2}.
 We observe a comparable tradeoff to that which may be gleaned from Section \ref{sec:experiments}: in particular, we are able to obtain complexity reduction by orders of magnitude with extremely little integral estimation error. This suggests the empirical validity of our compression-rule based on un-normalized importance weights operating in tandem with kernel smoothing.
%
%%%%%%%%%%%%%%%%%%%%%%%%%%%%%%%%%%%%%%%%%%%%%%%%%%%%%%%%%%%%%%%%%
\begin{figure*}[ht!]
	\centering
	\begin{subfigure}[b]{0.32\textwidth}
		\centering
		\includegraphics[width=\textwidth]{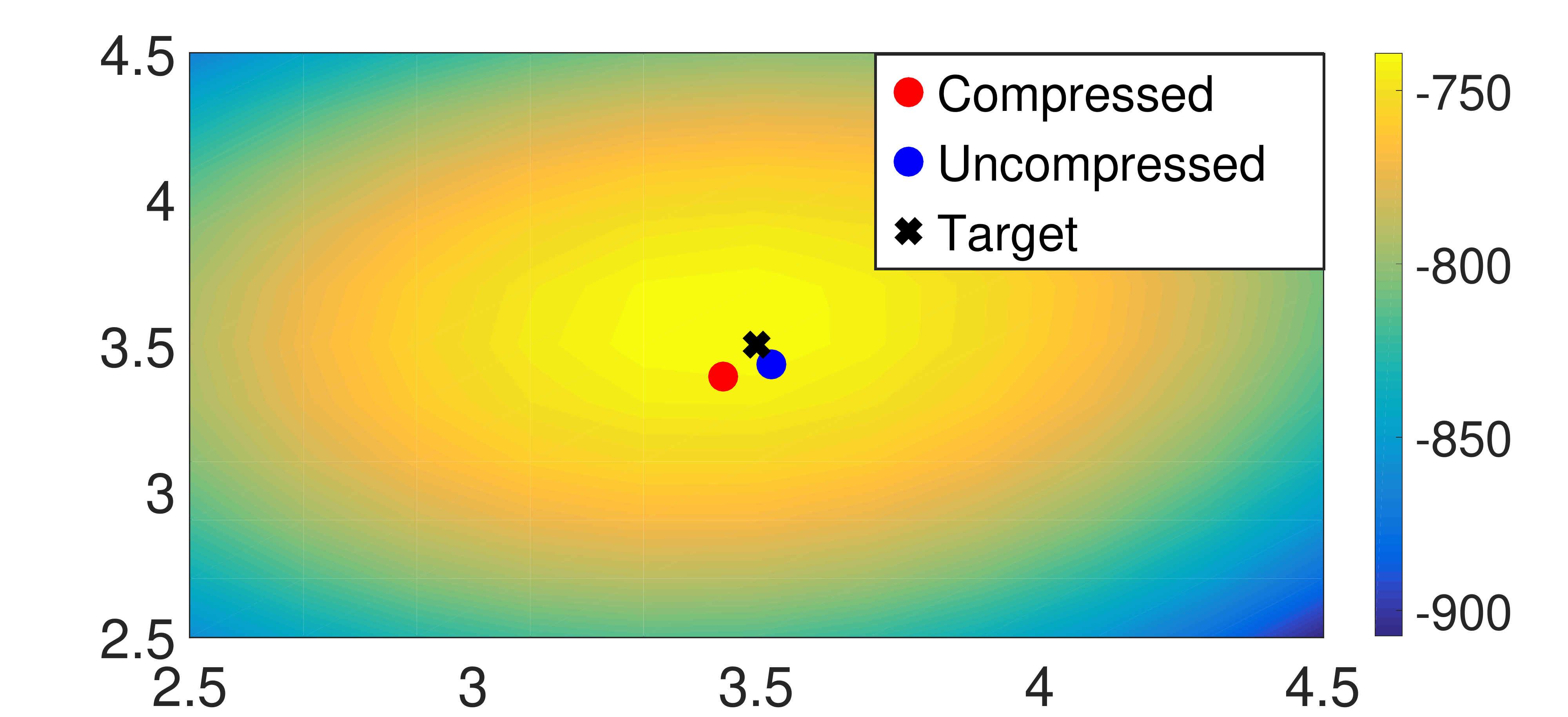}
		\caption{{\small Location estimate.}}    
		\label{fig:one}
	\end{subfigure}
	\begin{subfigure}[b]{0.32\textwidth}  
		\centering 
		\includegraphics[width=\textwidth]{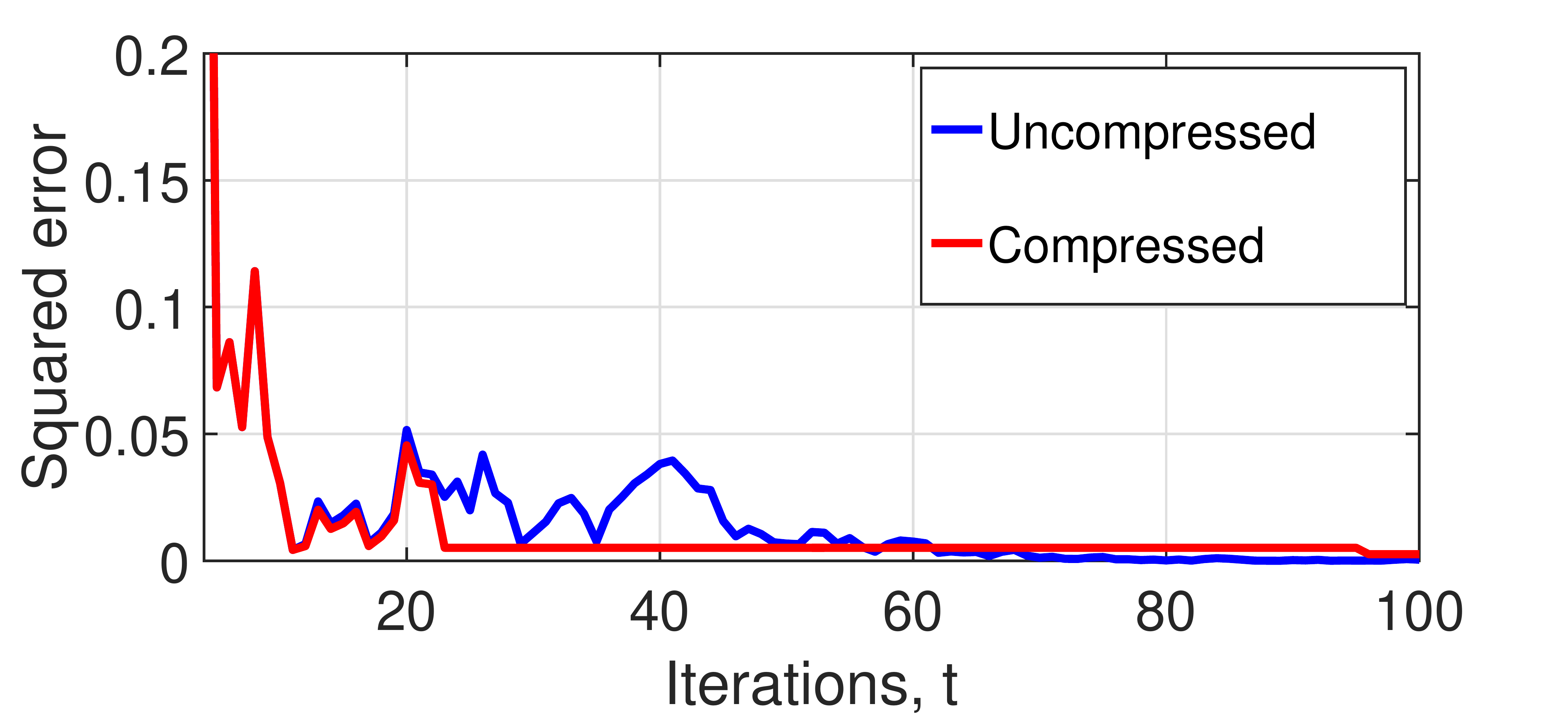}
		\caption[]%
		{{\small Squared error estimate.}}    
		\label{fig:two}
	\end{subfigure}
	\begin{subfigure}[b]{0.32\textwidth}   
		\centering 
		\includegraphics[width=\textwidth]{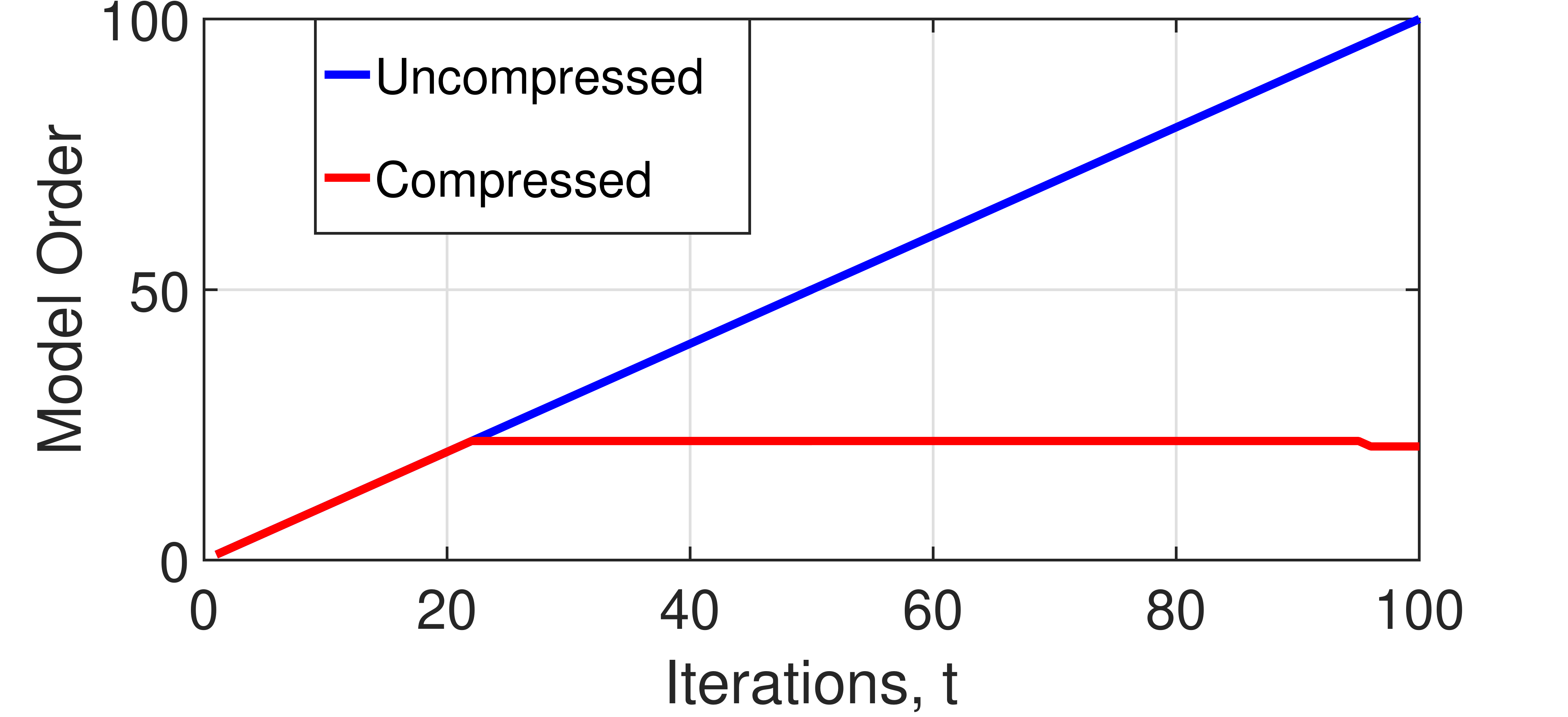}
		\caption[]%
		{{\small Model Order. }}    
		\label{fig:three}
	\end{subfigure}
	\caption[ The average and standard deviation of critical parameters ]
	{\small Simulation results for Alg. \ref{alg:compressed_online_kernel_importance_sampling} run with Gaussian kernel ($h=0.0001$) and compression budget $\epsilon=0.002$. Observe that the memory-reduction scheme (compressed) nearly preserves statistical consistency, while yielding a finite constant limiting model complexity, whereas the original uncompressed version (uncompressed) attains exact consistency but its memory grows linearly with particle index $t$. }
	\label{fig1}\vspace{-0mm}
\end{figure*}

\subsection{Source Localization}
In this section we present a sensor network localization experiment based on range measurements.The results illustrate the ability to succinctly  represent the unknown distribution of the source signal location, yielding a model that is both parsimonious and nearly consistent.
{Consider the problem of localizing a static target in two-dimensional space $\mathbb{R}^2$ with range measurements from the source collected in a wireless sensor network (WSN). Since the observation model is \emph{nonlinear}, the posterior distribution of the location of the target is intractable, and hence finding the least-squares estimator is not enough. This problem is well-studied in signal processing \cite{ali2009empirical} and robotics \cite{thrun2005probabilistic}.  Let $\x=[x,  y]^T$ denote the random unknown target location. We assume six sensors with locations $\{\bbh_i\}_{i=1}^6$ at locations $[1, -8]^T$, $[8, 10]^T$, $[-15, -17]^T$, $[-8, 1]^T$, $[10, 0]^T$, and $[0,10]^T$, respectively.  The true location of the target is at $[3.5, 3.5]^T$. The measurement at each sensor $i$ is related to the true target location $\x$ via the following nonlinear function of range
	%
	%\begin{align}
	%
	$y_{i,j}=-20\log(\|\x-\bbh_i\|)+\eta_i$
	%
	%\end{align}
	%
	for $i=1$ to $6$ and $j=1$ to $N$, where $N_i$  is the number of measurements collected by sensor $i$. Here, $\eta_i\sim\mathcal{N}(0,1)$ models the range estimation error.   For the experiment, we consider a Gaussian prior on the target location $\bbx$ with mean $[3.5,  3.5]^T$ and identity covariance matrix. We use  the actual target location as the mean for the Gaussian prior because we are interested in demonstrating that the proposed technique successfully balances particle growth and model bias. In practice, for a general possibly misspecified prior, we can appeal to advanced adaptive algorithms -- for example see \cite{elvira2017improving,bugallo2017adaptive}.}

Fig. \ref{fig1} shows the performance of the proposed algorithm compared against classical (uncompressed) normalized importance sampling. Fig. \ref{fig:one} shows that the final estimated value of the target location for compressed and uncompressed versions of the algorithms are close. We plot the squared error in Fig. \ref{fig:two} and both algorithms converge with close limiting estimates. Further, in Fig. \ref{fig:three} we observe that the model order for the compressed distribution settles to $21$, whereas the classical algorithm requires its number of particles in its importance distribution to grow unbounded. The memory-reduction comes at the cost of very little estimation error (Fig. \ref{fig:one}).   
\section{Conclusions}

We focused on Bayesian inference where one streams simulated Monte Carlo samples to approximate an unknown posterior via importance sampling. Doing so may consistently approximate any function of the posterior at the cost of infinite memory. Thus, we proposed Algorithm \ref{alg:compressed_online_kernel_importance_sampling} (CKIS) to approximate the posterior by a kernel density estimate (KDE) projected onto a nearby lower-dimensional subspace, which allows online compression as particles arrive in perpetuity. We established that the bias of CKIS depends on kernel bandwidth and compression budget, providing a tradeoff between statistical accuracy and memory. Experiments demonstrated that we attain memory-reduction by orders of magnitude with very little estimation error. This motivates future application to memory-efficient versions of Monte Carlo approaches to nonlinear signal processing problems such as localization, which has been eschewed due to its computational burden. %studied the compression bias that depends on an easily interpretable parameter $\epsilon$ that controls the the tradeoff between required memory and the compression loss. We have been able to prove that, under mild constraints, the compression loss vanishes when $\epsilon$ to zero (Theorem 1), and that the number of retained samples remains finite%

%%%%%%%%%%%%%%%%%%%%%%%%%%%%%%%%%%%%%%%%%%%%%%%%%%%%%%%%%%%%%%%%%%%%%%%%%%%%%%%%%%%%%%%%%%%%%%%%%%%%%%%%%%%%%%%%%%%%%%%%%%%%%%%%%%%%%%%%%%%%%%%%%%%%%%%S  E  C  T  I  O  N %%%%%%%%%%%%%%%%%%%%%%%%%%%%%%%%%%%%%%%%%%%%%%%%%%%%%%%%%%%%%%%%%%%%%%%%%%%%%%%%%%%%%%%%%%%%%%%%%%%%%%%%
%%%%%%%%%%%%%%%%%%%%%%%%%%%%%%%%%%%%%%%%%%%%%%%%%%%%%%%%%%%%%%%%%%
\appendices

\section{Proof of Consistency of Importance Sampling}\label{first_proof}
Here we state a result on the sample complexity and asymptotic consistency of IS estimators in terms of integral error. We increase the granularity of the proof found in the literature so that the modifications required for our results on compressed IS estimates are laid bare.
%\normalsize
%%%%%%%%%%%%%%%%%%%%%%%%%%%%%%%%%%%%%%%%%%%%%%%%%%%%%%%%%%%%%%%%%%%%%%%%%%%%%%%%%%%%%%%%%%%%%%%%%%%%%%%%%%%%%%%%%%%%%%%%%%%%%%%%%%%%%%%%%%%%%%%%%%%%%%% %%%%%%%%%%%%%%%%%%T  H  E  O  R  E  M%%%%%%%%%%%%%%%%%%%%%%%%%%%%%%%%%%%%%%%%%%%%%%%%%%%%%%%%%%%%%%%%%%%%%%%%%%%%%%%%%%%%%%
%%%%%%%%%%%%%%%%%%%%%%%%%%%%%%%%%%%%%%%%%%%%%%%%%%%%%%%%%%%%%%%%%%
 %%%%%%%%%%%%%%%%%%%%%%%%%%%%%%%%%%%%%%%%%%%%%%%%%%%%%%%%%%%%%%%%%%%%%%%%%%%%%%%%%%%%%%%%%%%%%%%%%%%%%%%%%%%%%%%%%%%%%%%%%%%%%%%%%%%%%%%%%%%%%%%%%%%%%%%%%%%%%%%%%%%%%%%%%%%%%%%%%%%%%%%%%%%%%%%%%%%
\begin{lemma}\label{first_theorem}
\cite{agapiou2017importance}[Theorem 2.1]  Suppose $\pi$, the proposal distribution is absolutely continuous w.r.t. $q$, the population posterior, and both are defined over $\ccalX$. Then define their Radon-Nikodyn derivative:
$
\frac{d q}{d \pi}(\bbx) := \frac{g(\bbx ) }{\int_{\ccalX} g(\bbx) \pi(d\bbx)} \; , \quad \rho := \frac{\pi(g^2)}{q(g^2)}
$
where $g$ is the unnormalized density of $q$ with respect to $\pi$. Moreover, $\rho$ is its second moment (``variance" of  unnormalized density).
Under Assumption \ref{as:first}\ref{as:testfunction1}, Alg. \ref{alg:importance_sampling} contracts to the true posterior as
\begin{align}\label{eq:posterior_contraction}
\sup_{|\phi| \leq 1}  | \mathbb{E}[I_N(\phi) - I(\phi) ] | \leq \frac{12}{N} \rho,\!\!\! 
\; \   \mathbb{E}\left[ ( I_N(\phi) - I(\phi) )^2\right]  \leq \frac{4}{N} \rho \; ,
\end{align}
and hence approaches exact consistency as $N \rightarrow \infty$.

\end{lemma}
%
%\footnote{Useful statements: $\rho \geq e^{D_{KL}(\mu \| \pi) }$ Theorem 4.19 in \cite{boucheron2013concentration}$D_{KL}(\mu \| \pi) \geq 2 d^2_H(\mu, \pi)$
%
%Therefore $\rho \geq e^{d^2_{H}(\mu , \pi) }$ which means that the distance from consistency scales exponentially with the square of the Hellinger metric } 
%
\begin{myproof}
This is a  more detailed proof than given in \cite{agapiou2017importance}[Theorem 2.1] develop for greater completeness and coherence. Let us denote the empirical random measure by $\pi_{}^N$ as
%
%\begin{align}\label{eq:occupancy_measure}
%
$\pi_{}^N:=\frac{1}{N}\sum\limits_{n=1}^{N}\delta_{\bbx{(n)}}$, and $\bbx{(n)}\sim \pi$,
%
%\end{align}
%
%\red{For the kernel density estimator, we write
%\begin{align}
%\pi_{}^N:=\frac{1}{N}\sum\limits_{n=1}^{N}\kappa_h({\bbx,\bbx{(n)}}), \ \ \ \bbx{(n)}\sim \pi.\end{align}
%where the kernel $\kappa_h({\bbx,\bbx{(n)}})$ is parametrized by the bandwidth parameter $h$. As $h\rightarrow 0$, the kernel function becomes the delta function (At least it holds for Gaussian kernel, we can start with the Gaussian itself right from the start and generalize later, if needed.)}
%
where $\pi_{}^N$ is the occupancy measure, which when weighted, yields the importance sampling empirical measure \eqref{eq:importance_sampling_empirical_measure}. Note that the integral approximation at $N$ is denoted by ${I}_N(\phi)$. With the above notation is hand, it holds that 
%\red{\begin{align}
%{I}_N(\phi)=\int\phi(\bbx)\pi_{}^N(\bbx)dx=\int\phi(\bbx)\frac{1}{N}\sum\limits_{n=1}^{N}\delta_{\bbx{(n)}}(\bbx)dx=\frac{1}{N}\sum\limits_{n=1}^{N}\phi({\bbx{(n)}})
%\end{align}}
%\red{Similiarily, for a kernel desity estimator, we can write
%\begin{align}
%{I}_N(\phi)=\int\phi(\bbx)\pi_{}^N(\bbx)dx=\int\phi(\bbx)\frac{1}{N}\sum\limits_{n=1}^{N}\kappa_h({\bbx,\bbx{(n)}})dx=\frac{1}{N}\sum\limits_{n=1}^{N}\phi({\bbx{(n)}})
%\end{align}where the last inequality holds from the reproducing property of the RKHS. Note that here we need the assumption that the function we are trying to estimate expected value of, belongs to RKHS, which is a valid one I guess.}
\begin{align}
\pi_{}^N(g)=\int\frac{1}{N}\sum\limits_{n=1}^{N}g(\bbx)\delta_{\bbx{(n)}}(\bbx
)dx=\frac{1}{N}\sum\limits_{n=1}^{N}g(\bbx{(n)}),
\end{align}
%\red{\begin{align}
%\pi_{}^N(g)=\int\frac{1}{N}\sum\limits_{n=1}^{N}g(\bbx)\delta_{\bbx{(n)}}(\bbx
%)dx=\frac{1}{N}\sum\limits_{n=1}^{N}g(x^{n}),
%\end{align}}
 and similarly 
\begin{align}
\pi_{}^N(\phi g)=&\int\frac{1}{N}\sum\limits_{n=1}^{N}\delta_{\bbx{(n)}}(\bbx
)\phi(\bbx
)g(\bbx)dx\nonumber \\
=&\frac{1}{N}\sum\limits_{n=1}^{N}\phi(\bbx{(n)})g(\bbx{(n)}).
\end{align}
%\red{\begin{align}
%\pi_{}^N(\phi g)=\int\frac{1}{N}\sum\limits_{n=1}^{N}\kappa(\bbx,{\bbx{(n)}})\phi(\bbx
%)g(\bbx)dx=\frac{1}{N}\sum\limits_{n=1}^{N}\phi(\bbx{(n)})g(\bbx{(n)}).
%\end{align}}
%
 From the above equalities, we can write the estimator bias as 
\begin{align}
{I}_N(\phi)-I(\phi)=&\frac{\pi_{}^N(\phi g)}{\pi_{}^N(g)}-I(\phi)
\\
=&\frac{\pi_{}^N(\phi g)}{\pi_{}^N(g)}-\left(I(\phi)\frac{\pi_{}^N(g)}{\pi_{}^N(g)}\right)
\\
=&\frac{1}{\pi_{}^N(g)}\left[\pi_{}^N(\phi g)-I(\phi)\pi_{}^N(g)\right]
\\
=&\frac{1}{\pi_{}^N(g)}\pi_{}^N\left((\phi-I(\phi))g\right).
\end{align}
Let us define $\bar{\phi}:=\phi-I(\phi)$ and note that 
\begin{align}\label{zero}
\pi(\bar{\phi}g)=0. \ \ \ 
\end{align}
Rewriting the bias, we get
\begin{align}\label{first}
{I}_N(\phi)-I(\phi)=&\frac{1}{\pi_{}^N(g)}\pi_{}^N\left(\bar{\phi}g\right)\nonumber
\\
=&\frac{1}{\pi_{}^N(g)}\left[\pi_{}^N\left(\bar{\phi}g\right)-\pi\left(\bar{\phi}g\right)\right],
\end{align}
where the second equality holds from \eqref{zero}. The first term in the bracket is an unbiased estimator for the second one, so that
\begin{align}\label{unbiased}
\mathbb{E}\left[\pi_{}^N\left(\bar{\phi}g\right)-\pi\left(\bar{\phi}g\right)\right]=0.
\end{align}
Taking the expectation on both sides  of \eqref{first}, we  get
\begin{align}\label{first1}
\mathbb{E}\left[{I}_N(\phi)-I(\phi)\right]=&\mathbb{E}\left[\frac{1}{\pi_{}^N(g)}
\left[\pi_{}^N\left(\bar{\phi}g\right)-\pi\left(\bar{\phi}g\right)\right]\right].
\end{align}
Since it equals zero, we can add the  expression in \eqref{unbiased} to the right hand side of \eqref{first1} to obtain
\begin{align}\label{first2}
\mathbb{E}\left[{I}_N(\phi)-I(\phi)\right]=&\mathbb{E}\left[\frac{1}{\pi_{}^N(g)}\left[\pi_{}^N\left(\bar{\phi}g\right)
-\pi\left(\bar{\phi}g\right)\right]\right]\nonumber\\
&\quad+\mathbb{E}\left[\pi_{}^N\left(\bar{\phi}g\right)-\pi\left(\bar{\phi}g\right)\right]
\\
&\hspace{-1cm}=\mathbb{E}\left[\frac{1}{\pi_{}^N(g)}\left[\pi_{}^N\left(\bar{\phi}g\right)-\pi\left(\bar{\phi}g\right)\right]\right]
\nonumber\\
&\hspace{-1cm}\quad+\mathbb{E}\left[\frac{1}{\pi(g)}\left(\pi_{}^N\left(\bar{\phi}g\right)-\pi\left(\bar{\phi}g\right)\right)\right].
\end{align}
Taking the expectation operator  outside, we get
\begin{align}
\mathbb{E}\left[{I}_N(\phi)\!-\!I(\phi)\right]=
&\mathbb{E}\!\left[\!\left(\!\frac{1}{\pi_{}^N(g)}\!-\!\frac{1}{\pi(g)}\!\right)\!\left(\pi_{}^N\left(\bar{\phi}g\right)-\pi\left(\bar{\phi}g\right)\right)\right]\nonumber\\
&\hspace{-2.8cm}=\mathbb{E}\left[\frac{1}{\pi_{}^N(g)\pi(g)}\left(\pi(g)-\pi_{}^N(g)\right)\left(\pi_{}^N\left(\bar{\phi}g\right)-\pi\left(\bar{\phi}g\right)\right)\right].
\end{align}
Next, we split the set of integration to \(A=\{2\pi^N_{MC}(g)>\pi(g)\}\) and its compliment using the property  
\[E[f(X)]=E[f(X)1_{A}(X)]+E[f(X)1_{A^c}(X)],\] where \(1_A\) the indicator function of the set \(A\) selecting $A=\{2\pi^N_{MC}(g)>\pi(g)\}$, which takes value 1 if \(x\in A\) and 0 if \(x\notin A\). We get
\begin{align}\label{intermediate}
|\mathbb{E}\left[{I}_N(\phi)-I(\phi)\right]|&\leq|\mathbb{E}\left[{I}_N(\phi)-I(\phi)\right]\boldsymbol{1}_{\{2\pi_{}^N(g)>\pi(g)\}}|  \nonumber\\
&\hspace{-1cm}+|\mathbb{E}\left[{I}_N(\phi)-I(\phi)\right]\boldsymbol{1}_{\{2\pi_{}^N(g)\leq\pi(g)\}}|.
\end{align}
Consider the second term of \eqref{intermediate}, and use the fact that $|\phi|\leq 1$, and so \(|\mu^N(\phi)|, |I(\phi)|\leq 1\) since they are mean values w.r.t. probability measures \(\mu^N,q\) respectively. Then we use \(E[1_A]=P(A)\) and obtain
\begin{align}
|\mathbb{E}\left[{I}_N(\phi)-I(\phi)\right]|\leq&|\mathbb{E}\left[{I}_N(\phi)
-I(\phi)\right]\boldsymbol{1}_{\{2\pi_{}^N(g)>\pi(g)\}}|\nonumber\\ 
&+2\mathbb{P}\left(2\pi_{}^N(g)\leq\pi(g)\right).\label{intermediate22}
\end{align}
The constant $2$ comes from the fact that $|{I}_N(\phi)-I(\phi)|\leq |{I}_N(\phi)|+|I(\phi)|\leq 2$. 
For the first term on the right hand side of \eqref{intermediate22}, from the set condition (, it holds that
\begin{align}
\frac{1}{\pi^N_{MC}(g)\pi(g)}<\frac{2}{\pi^2(g)},
\end{align}
which implies that
\begin{align}
|\mathbb{E}\left[{I}_N(\phi)\!-\!I(\phi)\right]|\leq&\frac{2}{\pi(g)^2}\mathbb{E}\left[|\pi(g)\!-\!\pi_{}^N(g)||\pi_{}^N(\bar{\phi}g)
\!-\!\pi\!\left(\bar{\phi}g\right)\!|\right]\nonumber\\
&+2\mathbb{P}\left(2\pi_{}^N(g)\leq\pi(g)\right).\label{intermediate12}
\end{align}
Finally, to upper bound the first term on the right hand side of \eqref{intermediate12}, we first bound the expectation  using Cauchy-Schwartz 
\begin{align}
E[|\pi(g)-&\pi_{}^N(g)||\pi^N(\bar{\phi}g)-\pi(\bar{\phi}g)|]
\nonumber\\
&\leq E[(\pi(g)-\pi_{}^N(g))^2]^{\frac12}E[(\pi^N(\bar{\phi}g)-\pi(g))^2]^\frac12\label{cauchy1}
\end{align}
The first expectation on the right hand side of \eqref{cauchy1} is bounded as follows: 
by definition of \(\pi^N\) we have for \(x_n\sim \pi\) independent that
 \begin{align}
 E[(\pi(g)-\pi_{}^N(g))^2]=& E[(\pi(g)-\frac1N\sum_{n=1}^Ng(\bbx{(n)})^2)]\nonumber \\
 &\hspace{-1cm}=\frac1{N^2}E[(\sum_{n=1}^Ng(\bbx{(n)})-N\pi(g))^2],
 \end{align} 
 which since \(E[g(\bbx{(n)})]=\pi(g)\) 
 and by independence of the \(\bbx{(n)}\) is equal to \[=\frac1{N^2} Var(\sum_{n=1}^Ng(\bbx{(n)}))=\frac1{N^2} \sum_{n=1}^NVar(g(\bbx{(n)})),\] and since \(\bbx{(n)}\) is identically distributed, \(\bbx{(n)}\sim\pi\), 
 this is equal to \[\frac{N}{N^2}Var_{u\sim \pi}(g)=\frac1N(\pi(g^2)-\pi(g)^2)\leq \frac1N\pi(g^2).\]
The second expectation on the right hand side of \eqref{cauchy1} is bounded in a similar way along with the fact that \(|\phi|\leq 1\) so that \(|\bar{\phi}|\leq 2\). Then we utilize these upper bounds on the right hand side of \eqref{intermediate} to obtain
\begin{align}
|\mathbb{E}\left[{I}_N(\phi)\!-\!I(\phi)\right]|\leq&\frac{2}{\pi(g)^2}\mathbb{E}\left[|\pi(g)\!-\!\pi_{}^N(g)||\pi_{}^N(\bar{\phi}g)
\!-\!\pi\left(\!\bar{\phi}g\!\right)\!|\right]\nonumber\\
&+2\mathbb{P}\left(2\pi_{}^N(g)\leq\pi(g)\right)\label{intermediate1}\\
\leq&\frac{2}{\pi(g)^2}\frac{1}{\sqrt{N}}\pi(g^2)^{1/2}\frac{2}{\sqrt{N}\pi(g^2)^{1/2}}]\nonumber\\
&+2\mathbb{P}\left(2\pi_{}^N(g)\leq\pi(g)\right)\label{intermediate2}
\end{align}
where the inequalities follow from the fact that the test function is bounded $|\phi|$. Next, note that
\begin{align}\label{prob}
\mathbb{P}\left(2\pi_{}^N(g)\leq\pi(g)\right)=&\mathbb{P}\left(2\left(\pi_{}^N(g)-\pi(g)\right)\leq-\pi(g)\right)\nonumber\\
\leq
& \mathbb{P}\left(2|\pi_{}^N(g)-\pi(g)|\geq\pi(g)\right),
\end{align} 
where the first equality  is obtained by subtracting $-2\pi(g)$ from both sides inside the bracket.
Next, we use the Markov inequality, given by
%
%\begin{align}
%
$\mathbb{P}\left(X\geq a\right)\leq \frac{\mathbb{E}(X)}{a}$.
%
%\end{align}
%
Utilizing this, we can write
\begin{align*}
\mathbb{P}\left(2|\pi_{}^N(g)\!-\!\pi(g)|\geq\pi(g)\right)
\!\leq \frac{2\mathbb{E}\left[|\pi_{}^N(g)\!-\!\pi(g)|\right]}{\pi(g)}
\!\!\leq\!\! \frac{4}{N}\frac{\pi(g^2)}{\pi(g)^2}.
\end{align*}
This implies that
\begin{align}
\mathbb{P}\left(2\pi_{}^N(g)\leq\pi(g)\right)\leq \frac{4}{N}\frac{\pi(g^2)}{\pi(g)^2}.\label{inter2}
\end{align}
Finally, using  the upper bound in \eqref{inter2} in \eqref{intermediate2}, we obtain
\begin{align}
\sup_{|\phi|\leq 1}|\mathbb{E}\left[{I}_N(\phi)-I(\phi)\right]|\leq&\frac{12}{N}\frac{\pi(g^2)}{\pi(g)^2}
\end{align}
which proves the result.% $\hfill \blacksquare$
\end{myproof}

%%%%%%%%%%%%%%%%%%%%%%%%%%%%%%%%%%%%%%%%%%%%%%%%%%%%%%%%%%%%%%%%%%%%%%%%%%%%%%%%%%%%%%%%%%%%%%%%%%%%%%%%%%%%%%%%%%%%%%%%%%%%%%%%%%%%%%%%%%%%%%%%%%%%%% %%%%%%%%%%%%%%%%%%P  R  O  O  F%%%%%%%%%%%%%%%%%%%%%%%%%%%%%%%%%%%%%%%%%%%%%%%%%%%%%%%%%%%%%%%%%%%%%%%%%%%%%%%%%%%%%%
%%%%%%%%%%%%%%%%%%%%%%%%%%%%%%%%%%%%%%%%%%%%%%%%%%%%%%%%%%%%%%%%%%
%\section{Proof of Theorem \ref{thm:main_result}}\label{second_proof}

%
\section{Proof of Theorem \ref{thm:model_order}}\label{third_proof}
We begin by presenting a lemma which allows us to relate the stopping criterion of our sparsification procedure to a Hilbert subspace distance.

%\red{discussion of lemma}
%%%%%%%%%%%%%%%%%%%%%%%%%%%%%%%%%%%%%%%%%%%%%%%%%%%%%%%%%%%%%%%
%%% L  E  M  M  A  %%%%%%%%%%%%%%%%%
%%%%%%%%%%%%%%%%%%%%%%%%%%%%%%%%%%%%%%%%%%%%%%%%%%%%%%%%%%%%%%%
%
\begin{lemma}\label{lemma_subspace_dist}
Define the distance of an arbitrary feature vector $\bbx$ evaluated by the feature transformation $\psi(\bbx):=\kappa(\bbx,\cdot)$ to $\ccalH_{\bbD}=\text{span}\{\psi(\bbd_n) \}_{n=1}^M$, the subspace of the real space spanned by a dictionary $\bbD$ of size $M$,  as
\begin{align}\label{eq:hilbert_subspace_dist}
\text{dist}( \psi(\bbx) , \ccalH_{\bbD}) 
= \min_{y\in\ccalH_{\bbD}} | \psi(\bbx) - \bbv^T \bbphi_{\bbD} | \; .
\end{align}
This set distance simplifies to the following least-squares projection when $\bbD \in \reals^{p\times M}$ is fixed
\begin{align}\label{eq:hilbert_subspace_dist_ls}
\text{dist}( \psi(\bbx) , \ccalH_{\bbD}) 
=& \Big|  \psi(\bbx) 
  - \psi(\bbx)\bbpsi_{\bbD}^T\bbK_{\bbD, \bbD}^{-1} 
   \bbpsi_{\bbD} \Big|
   .
\end{align}
\end{lemma}

%%%%%%%%%%%%%%%%%%%%%%%%%%%%%%%%%%%%%%%%%
%%%%%%%%%%%%%%%%%%%%%%%%%%%%%%%%%%%%%%%%%
%%%%%%%%    P   R   O   O   F    %%%%%%%%%%%%%%%%%%%%%%% 
%%%%%%%%%%%%%%%%%%%%%%%%%%%%%%%%%%%%%%%%%
%%%%%%%%%%%%%%%%%%%%%%%%%%%%%%%%%%%%%%%%%
%\begin{proof} See Appendix 

%We now present our main result, which says that the model order of $f_n$ remains bounded for all $t$, and depends on the conditioning of the problem setting as well as the feature space.

%\subsection{Proof of Lemma \ref{lemma_subspace_dist}}\label{apx_lemma_subspace_dist}
%%%%%%%%%%%%%%%%%%%%%%%%%%%%%%%%%%%%%%%%%
%%%%%%%%%%%%%%%%%%%%%%%%%%%%%%%%%%%%%%%%%
%%%%%%%%    P   R   O   O   F    %%%%%%%%%%%%%%%%%%%%%%% 
%%%%%%%%%%%%%%%%%%%%%%%%%%%%%%%%%%%%%%%%%
%%%%%%%%%%%%%%%%%%%%%%%%%%%%%%%%%%%%%%%%%
\begin{myproof}
The distance to the subspace $\ccalH_{\bbD}$ is defined as
\begin{align}\label{eq:subspace_dist}
\text{dist}( \psi(\bbx) , \ccalH_{\bbD_n}) 
=& \min_{y\in\ccalH_{\bbD}} | \psi(\bbx) - \bbv^T \bbpsi_{\bbD}|
\nonumber\\=& \min_{\bbv\in \reals^{M}} | \psi(\bbx) - \bbv^T \bbpsi_{\bbD} | \; ,
\end{align}
where the first equality comes from the fact that the dictionary $\bbD$ is fixed, so $\bbv\in \reals^M$ is the only free parameter. Now plug in the minimizing weight vector $\tbv^\star=\psi(\bbx_n)\bbK_{\bbD_n, \bbD_n}^{-1}\bbpsi_{\bbD_n}$ into \eqref{eq:subspace_dist} which is obtained in an analogous manner to the logic which yields \eqref{eq:proximal_hilbert_representer} - \eqref{eq:hatparam_update}. Doing so simplifies \eqref{eq:subspace_dist} to the following
\begin{align}\label{eq:subspace_dist2}
\text{dist}(\psi(\bbx_n) , \ccalH_{\bbD_n}) =&  \Big|  \psi(\bbx_n) 
  - \psi(\bbx_n)[\bbK_{\bbD_n, \bbD_n}^{-1} \bbpsi_{\bbD_n}]^T
   \bbpsi_{\bbD_n} \Big|\nonumber \;\\
   &\hspace{-10mm}=  \Big|  \psi(\bbx_n) 
     - \psi(\bbx_n)\bbpsi_{\bbD_n}^T\bbK_{\bbD_n, \bbD_n}^{-1}
      \bbpsi_{\bbD_n} \Big| .
\end{align}
\end{myproof}

Next, we establish that the model order is finite.

\begin{myproof}
Consider the model order of the kernel mean embedding $\beta_{n}$ and $\beta_{n-1}$ generated by Algorithm \ref{alg:compressed_online_kernel_importance_sampling} and denoted by $M_{n}$ and $M_{n-1}$, respectively, at two arbitrary subsequent instances $n$ and $n-1$.
% Assume a constant algorithm step-size $\eta$ has been chosen such that $\eta<1/\lambda$ and the approximation budget $\eps$ satisfies $\eps=K \eta^{3/2}$ for some positive scalar $K>0$.
%
Suppose the model order of the estimate $\beta_{n}$ is less than or equal to that of $\beta_{n-1}$, i.e. $M_{n} \leq M_{n-1}$. This relation holds when the stopping criterion of MMD-OMP ( defined in Algorithm \ref{alg:compressed_online_kernel_importance_sampling}), stated as $\min_{j=1,\dots,{M_{n-1} + 1}} \gamma_j > \eps$, \emph{is not} satisfied for the updated dictionary matrix with the newest sample point $\bbx{(n)}$ appended: $\tbD_{n} = [\bbD_{n-1} ; \bbx{(n)} ]$ [cf. \eqref{eq:importance_sampling_stack2}], which is of size $M_{n-1} + 1$. 
Thus, the negation of the termination condition of MMD-OMP in Algorithm \ref{alg:compressed_online_kernel_importance_sampling} must hold for this case, stated as
\begin{align}\label{eq:komp_no_terminate}
\min_{j=1,\dots,{M_{n-1} + 1}} \gamma_j \leq \eps \; .
\end{align}
Observe that the left-hand side of \eqref{eq:komp_no_terminate} lower bounds the approximation error $\gamma_{M_{n-1} + 1}$ for removing the most recent sample $\bbx{(n)}$ due to the minimization over $j$, that is, $\min_{j=1,\dots,{M_{n-1} + 1}} \gamma_j \leq\gamma_{M_{n-1} + 1} $. Consequently, if $\gamma_{M_{n-1} + 1} \leq \eps$, then \eqref{eq:komp_no_terminate} holds and the model order does not grow. Thus it suffices to consider $\gamma_{M_{n-1} + 1}$.
The definition of $\gamma_{M_{n-1} + 1}$ with the substitution of $\beta_{n}$ in \eqref{eq:komp_no_terminate} allows us to write
\begin{align}\label{eq:min_gamma_expand}
\gamma_{M_{n-1}+1}
			&\!\!=\!\!\!\!\!\min_{\bbu\in\reals^{{M_{n-1}}}} \!\!\Big| {\beta_{n-1}} +  g(\bbx{(n)}) \kappa_{\bbx{(n)}}(\bbx)  -\hspace{-0.8cm}\sum_{k \in \ccalI \setminus \{M_{n-1} + 1\}} \hspace{-0.8cm}u_k \kappa_{\bbd_k}(\bbx) \Big| \nonumber\\
			&=\min_{\bbu\in\reals^{{M_{n-1}}}} \Big|  \sum_{k \in \ccalI \setminus \{M_{n-1} + 1\}} \!\!\!  \!\!\!g(\bbx(k)) \kappa_{\bbd_{k}}(\bbx)\\
			&\qquad\qquad\quad\!\!+ g(\bbx{(n)}) \kappa_{\bbx{(n)}}(\x)   -\hspace{-0.8cm} \sum_{k \in \ccalI \setminus \{M_{n-1} + 1\}} \hspace{-0.8cm}u_k \kappa_{\bbd_{k}}(\bbx)\Big |\; \nonumber
			, 
\end{align}  
where we denote $\kappa_{\bbx{(n)}}(\bbx)=\kappa_{(\bbx{(n)},\cdot)}$ and the $k^\text{th}$ column of $\bbD_n$ as $\bbd_k$. The minimal error is achieved by considering the square of the expression inside the minimization and expanding terms to obtain
\begin{align}\label{eq:error_expansion}
\Big| &\!\!\! \!\!\!\sum_{k \in \ccalI \setminus \{M_{n-1} + 1\}} \!\!\!  \!\!\!g(\bbx(k)) \kappa_{\bbd_{k}}(\bbx)+  g(\bbx{(n)}) \kappa_{\bbx{(n)}}(\x)   - \hspace{-1cm} \sum_{k \in \ccalI \setminus \{M_{n-1} + 1\}} \hspace{-0.8cm}u_k \kappa_{\bbd_{k}}(\bbx)\Big|^2  \nonumber\\
&=\Big|  \bbg^T\bbkappa_{\bbD_n}(\bbx) + g(\bbx{(n)}) \kappa_{\bbx{(n)}}(\x)  -  \bbu^T \bbkappa_{\bbD_n}(\bbx) \Big|^2 \nonumber\\
&=  \bbg^T \bbK_{\bbD_n, \bbD_n} \bbg 
+  g(\bbx{(n)})^2 
+ \bbu^T \bbK_{\bbD_n, \bbD_n} \bbu  \nonumber \\
&\quad\!+\! 2 g(\bbx{(n)}) \bbw^T\! \bbkappa_{\bbD_n}(\bbx{(n)})\! 
\!-\! 2  g(\bbx{(n)}) \bbu^T\!\bbkappa_{\bbD_n}(\bbx{(n)})\nonumber \\
&\qquad-2  \bbw^T\! \bbK_{\bbD_n, \bbD_n}\! \bbu  . 
\end{align} 
 To obtain the minimum, we compute the stationary solution of \eqref{eq:error_expansion} with respect to $\bbu \in \reals^{M_{n-1}}$ and solve for the minimizing $\tbu^\star$, which in a manner similar to the logic in \eqref{eq:proximal_hilbert_representer} - \eqref{eq:hatparam_update}, is given as
%
% \begin{align}\label{eq:minimal_weights}
%
$ \tbu^\star =  {[ g(\bbx{(n)})\bbK_{\bbD_n, \bbD_n}^{-1}\bbkappa_{\bbD_n}(\bbx{(n)})
 +   \bbg ]} \; .$
%\end{align}
%
Plug $\tbu^\star$ into the expression in \eqref{eq:min_gamma_expand} and,  using the short-hand notation $\sum_k u_k \kappa_{\bbd_{k}}(\bbx)= \bbu^T \bbkappa_{\bbD_n}(\bbx)$. Simplifies \eqref{eq:min_gamma_expand} to
\begin{align}\label{eq:min_gamma_optimal_weights}
&\Big| \bbg^T\bbkappa_{\bbD_n}(\bbx) + g(\bbx{(n)}) \kappa_{\bbx{(n)}}(\x)  -  \bbu^T \bbkappa_{\bbD_n}(\bbx) \Big| \nonumber\\
& =\! \Big| \bbg^T\bbkappa_{\bbD_n}(\bbx) + g(\bbx{(n)}) \kappa_{\bbx{(n)}}(\x) \! \nonumber \\
& -  {[ g(\bbx{(n)})\bbK_{\bbD_n, \bbD_n}^{-1} \bbkappa_{\bbD_n}(\bbx{(n)})\! \!
 + \!\!  \bbg ]}^T 
\!\!\bbkappa_{\bbD_n}(\bbx) \! \Big|  \nonumber\\
& =\! \Big| g(\bbx{(n)}) \kappa_{\bbx{(n)}}(\x)\! -\! {[ g(\bbx{(n)})  \bbK_{\bbD_n, \bbD_n}^{-1}\bbkappa_{\bbD_n}(\bbx{(n)})
 ]}^T \!\!
\bbkappa_{\bbD_n}(\bbx)\Big|  \nonumber\\
& =\! g(\bbx{(n)}) \Big| \kappa_{\bbx{(n)}}(\x)-  \bbkappa_{\bbD_n}(\bbx^n)^T\bbK_{\bbD_n, \bbD_n}^{-1} \bbkappa_{\bbD_n}(\bbx)\Big|
\end{align}
Notice that the right-hand side of \eqref{eq:min_gamma_optimal_weights} may be identified as the distance to the subspace $\ccalH_{\bbD_n}$ in \eqref{eq:subspace_dist2} defined in Lemma \ref{lemma_subspace_dist} scaled by a factor of {$g(\bbx{(n)})$}. We may upper-bound the right-hand side of \eqref{eq:min_gamma_optimal_weights} as
\begin{align}\label{eq:min_gamma_optimal_weights3}
   g(\bbx{(n)}) &\Big| \kappa_{\bbx{(n)}}(\x)-  \bbkappa_{\bbD_n}(\bbx{(n)})^T \bbK_{\bbD_n, \bbD_n}^{-1}\bbkappa_{\bbD_n}(\bbx)\Big| \nonumber\\
   &=   g(\bbx{(n)})\text{dist}(\kappa_{\bbx{(n)}}(\bbx) , \ccalH_{\bbD_n})
\end{align}
where we have applied \eqref{eq:hilbert_subspace_dist_ls} regarding the definition of the subspace distance on the right-hand side of \eqref{eq:min_gamma_optimal_weights3} to replace the absolute value term. Now, when the {MMD-OMP} stopping criterion is violated, i.e., \eqref{eq:komp_no_terminate} holds, this implies $\gamma_{M_{n-1} + 1} \leq \eps$. Therefore, the right-hand side of \eqref{eq:min_gamma_optimal_weights3} is upper-bounded by $\epsilon$, and we can write
\begin{align}\label{eq:min_gamma_optimal_weights5}
 g(\bbx{(n)})\text{dist}(\kappa_{\bbx{(n)}}(\bbx) , \ccalH_{\bbD_n})\leq \epsilon.
\end{align}
 After rearranging the terms in \eqref{eq:min_gamma_optimal_weights5}, we can write
\begin{align}\label{eq:min_gamma_optimal_weights4}
  \text{dist}(\kappa_{\bbx{(n)}}(\bbx),\ccalH_{\bbD_n}) \leq \frac{\epsilon}{ g(\bbx{(n)})} \;,
\end{align}
where we have divided both sides by {$g(\bbx{(n)})$}. 
%
%Using this identification, we transform the sufficient condition for the stopping criterion of KOMP to be violated, stated as $\gamma_{M_{n-1} + 1} \leq \eps$, into a criterion on $\text{dist}(\kappa_{\bbx{(n)}}(\bbx),\ccalH_{\bbD_n})$, the subspace distance of $\kappa_{\bbx{(n)}}(\bbx)$ to the span of kernel evaluations of the current dictionary $\ccalH_{\bbD_n}$.
%%
%%
%%\red{more explanation here}
%Substituting the definition \eqref{eq:subspace_dist2} into $\gamma_{M_{n-1} + 1} \leq \eps$ and dividing both sides by $\eta |\ell_n'({f}_n(\bbx{(n)})) |$ yields
%%
%\begin{align}\label{eq:newest_gamma_rescaled}
%\text{dist}(\kappa_{\bbx{(n)}}(\bbx),\ccalH_{\bbD_n}) \leq \frac{\epsilon}{\eta|\ell'(f_n(\bbx{(n)}),\bby_n)|} \; .
%\end{align}
%%
%%Now use the approximation budget selection in terms of the learning rate $\eta$ as $\epsilon=K\eta^{3/2}$. Furthermore, the $C$-Lipschitz continuity of $\ell$ [cf. \eqref{eq:lipschitz}] in Assumption \ref{as:2} allows us to bound the instantaneous gradient by this same constant. Inverting this expression yields $1/|\ell'(f_n(\bbx{(n)}),\bby_n)| \geq 1/C$. Substituting in this lower bound and selection of $\eps$, we obtain that if
%Now, divide the above expression by $C \eta$ to write
%%%
%\begin{align}\label{eq:newest_gamma_rescaled2}
%%
%\text{dist}(\kappa_{\bbx{(n)}}(\bbx),\ccalH_{\bbD_n}) \leq \frac{K\sqrt{\eta}}{C}
%%
%\end{align}
%%
Observe that if \eqref{eq:min_gamma_optimal_weights4} holds, then $\gamma_{M_{n}} \leq \eps$ holds, but since $\gamma_{M_{n}} \geq \min_{j} \gamma_j $, we may conclude that \eqref{eq:komp_no_terminate} is satisfied. Consequently the model order at the subsequent step does not grow which means that $M_{n} \leq M_{n-1}$ whenever \eqref{eq:min_gamma_optimal_weights4} is valid. 

Now, let's take the contrapositive of the preceding expressions to observe that growth in the model order ($M_{n} = M_{n-1} + 1$) implies that the condition
\begin{align}
\label{eq:min_gamma2}
\text{dist}(\kappa_{\bbx{(n)}}(\bbx),\ccalH_{\bbD_n}) > \frac{\epsilon}{ g(\bbx{(n)})} 
\end{align} %\vspace{-5cm}
holds.  Therefore, each time a new point is added to the model, the corresponding map $\kappa_{\x}(\bbx{(n)})$ is guaranteed to be at least a distance of {$\frac{\epsilon}{g(\bbx{(n)})} $} from every other feature map in the current model. 
In canonical works such as \cite{bengtsson2008curse,bickel2008sharp}, the largest self-normalized importance weight is shown to be bounded by a constant. Under the additional hypothesis that the \emph{un-normalized} importance weight is bounded by some constant $W$, then we have via \eqref{eq:min_gamma2}
%
%{\begin{align}
%\label{eq:min_gamma3}
$\text{dist}(\kappa_{\bbx{(n)}}(\bbx),\ccalH_{\bbD_n}) > \frac{\epsilon}{W} .$
%\end{align}}
%
 Therefore, For a fixed compression budget $\epsilon$, {the MMD-OMP} stopping criterion is violated for the newest point whenever distinct dictionary points $\bbd_k$ and $\bbd_j$  for $j,k\in\{1,\dots,M_{n-1}\}$, satisfy the condition $\text{dist}(\kappa_{\x}(\bbd_j),\kappa_{\bbd_{k}}(\bbx)) > \frac{\epsilon}{W}$. Next, we follow a similar argument as provided in the proof of Theorem 3.1 in \cite{1315946}.   Since $\ccalX$ is compact and $\kappa_{\x}$ is continuous, the range $\kappa_{\x}(\ccalX) $ of the feature space $\ccalX$ is compact.   Therefore, the minimum number of balls (covering number) of radius $\kappa$ (here, $\kappa = \frac{\epsilon}{W}$) needed to cover the set $\kappa_{\x}(\ccalX)$ is finite (see, e.g., \cite{anthony2009neural}) for a finite compression budget $\epsilon$. The finiteness of the covering number implies that the number of elements in the dictionary $M_N$ will be finite  and using \cite[Proposition 2.2]{1315946}, we can characterize the number of elements in the dictionary as
 %
%\begin{align}
$1 \leq M_{N}\leq C\left(\frac{W}{\eps}\right)^{2p}$,
%\end{align}
 %
where $C$ is a constant depending upon the space $\mathcal{X}$.
\end{myproof}\vspace{-0mm}

%\small
\bibliographystyle{IEEEtran}
\bibliography{bibliography}
%\normalsize

%\clearpage
%\newpage

   \end{document}